\documentclass[a4paper]{amsart}
\usepackage{amsmath, amsfonts,amsthm,amssymb,amscd, verbatim,graphicx,color,multirow,booktabs, caption, mathdots,bm,chngcntr,adjustbox,longtable,mathtools,microtype,booktabs, amsthm}
\usepackage{pgfplots}
\usepackage{xcolor}
\usepackage[top = 78pt, bottom = 72pt, inner=72pt, outer=72pt]{geometry}

\definecolor{Teal}{HTML}{23373b}
\definecolor{Purple}{HTML}{606BA9}
\definecolor{Orange}{HTML}{DF7C2A}
\definecolor{Red}{HTML}{BD0808}
\definecolor{Lightpurple}{HTML}{7180AC}
\definecolor{BurntOrange}{HTML}{DD6E42}
\definecolor{Lighterblue}{HTML}{64CCE9}
\definecolor{lavander}{HTML}{B56BEE}
\definecolor{Seagreen}{HTML}{39932F}

\usepackage{hyperref}
\hypersetup{colorlinks=true, linkcolor=BurntOrange, citecolor=BurntOrange}
\usepackage[capitalise,noabbrev]{cleveref}
\usepackage[capitalise,noabbrev]{cleveref}
\usepackage[shortlabels]{enumitem}

\newtheorem{theorem}{Theorem}[]
\newtheorem{lemma}[theorem]{Lemma}
\newtheorem{cor}[theorem]{Corollary}

\newtheorem{condition}[theorem]{Condition}
	
\theoremstyle{definition}
\newtheorem*{de}{Definition}
\newtheorem{remark}[theorem]{Remark}

\newtheorem{proposition}[theorem]{Proposition}

\newcommand{\V}{\mathrm{V}}
\newcommand{\E}{\mathrm{E}}

\newcommand{\Aut}{\mathrm{Aut}}
\newcommand{\NN}{\mathbb{N}}
\newcommand{\ZZ}{\mathbb{Z}}
\renewcommand{\l}{\ell}
\newcommand{\A}{\mathbf A}
\renewcommand{\O}{\mathcal{O}}
\renewcommand{\r}{\mathbf{r}}
\DeclareMathOperator{\Pet}{Pet}
\DeclareMathOperator{\Dih}{Dih}
\newcommand*{\lmulti  }{\{\mskip-5mu\{}
\newcommand*{\rmulti}{\}\mskip-5mu\}}

\newcommand{\C}{\mathcal{C}}

\usepackage{tikz}
\usetikzlibrary{graphs}
\usepackage{adjustbox}
\usepackage{subcaption}
\usetikzlibrary{graphs,graphs.standard, calc}
\tikzset{normaledge/.style={inner sep=2pt,fill,outer sep=0,circle}}
\tikzset{rededge/.style={inner sep=2pt,fill,outer sep=0,circle}}

\pgfplotsset{compat=1.18}

\tikzstyle{krog}=[draw=black, minimum size=0.65cm, shape=circle]
\tikzstyle{pika}=[fill=black, draw=none, shape=circle, scale=0.1]
\tikzstyle{none}=[fill=none, draw=none, shape=circle, scale=0.1]

\tikzstyle{povezava}=[-, fill=none]
\tikzstyle{rdeca povezava}=[-, draw=Red]
\tikzstyle{crtkana povezava}=[-, draw=black, dashed, dash pattern=on 4mm off 2mm]

\pgfdeclarelayer{nodelayer}
\pgfdeclarelayer{edgelayer}
\pgfsetlayers{edgelayer,main,nodelayer}

\newsavebox\mybox
\newenvironment{resizedtikzpicture}[1]{%
  \def\mywidth{#1}%
  \begin{lrbox}{\mybox}%
  \begin{tikzpicture}
}{%
  \end{tikzpicture}%
  \end{lrbox}%
  \resizebox{\mywidth}{!}{\usebox\mybox}%
}

\title{Cubic vertex-transitive graphs of girth seven}

\usepackage[%
sorting=nyt, 
 style=numeric,
 url=false,
maxbibnames=9,%
 doi=false,
 isbn=false,
sortcites = true,
citestyle=numeric-comp, 
]{biblatex}

\author[M. ~Lekše ]{Maruša Lekše}
\address{Institute of Mathematics, Physics, and Mechanics, Jadranska ulica 19, 1000 Ljubljana, Slovenia. Also affiliated with: Faculty of Mathematics and Physics, University of Ljubljana, Jadranska ulica 21, 1000 Ljubljana, Slovenia.} 
\email{marusa.lekse@imfm.si}

\author[M.~Toledo]{Micael Toledo}
\address{Faculty of Mathematics and Physics, University of Ljubljana, Jadranska ulica 21, 1000 Ljubljana, Slovenia. Also affiliated with: Institute of Mathematics, Physics, and Mechanics, Jadranska ulica 19, 1000 Ljubljana, Slovenia.} 
\email{micael.toledo@fmf.uni-lj.si}

\begin{document}

\maketitle

\begin{abstract}
In this paper we classify cubic vertex-transitive graphs of girth $7$, based on their signature. Such a graph is either a truncation of an arc-transitive dihedral scheme on a $7$-regular graph, the skeleton of a rotary map of type $\{7,3\}$, a member of an infinite family of Cayley graphs, or is one of the of the generalised Petersen graphs $\text{Pet}(13,5)$, $\text{Pet}(15,4)$, $\text{Pet}(17,4)$ or the Coxeter graph. We show that for a cubic vertex-transitive graphs $\Gamma$ of girth $7$, if every edge of $\Gamma$ is contained in the same number of $7$-cycles, then $\Gamma$ is also arc-transitive.
\end{abstract}

\section{Introduction}

    The presence of symmetry in graphs constrains their structure and can be exploited to facilitate their study, often allowing for full classification results. Cubic (trivalent) vertex-transitive graphs are, arguably, the most studied class of highly symmetric graphs and have attracted the attention of notable names like Tutte, Frucht, Foster and Coxeter, to name a few (see for instance \cite{FosterrCensus,Tutte_1947,LORIMER198337,zerosym} for a very small sample of the classic work in this area). Currently, there exist several censuses of cubic vertex-transitive graphs \cite{Potocnik+Spiga+Verret+2013,PotocnikCayley5000,ConderCensus} and an even larger number of classification results of cubic vertex-transitive graphs with a variety of prescribed conditions; be it a geometric property \cite{MARK2017}, the presence of a semiregular automorphism \cite{PISANSKI2007567,POTOCNIKTOLEDO2020}, the presence of a subgroup of its automorphism group with a given structure or type of action \cite{POTOCNIKTOLEDO2021,ZHOU2014679,NonCayely}, the order of the graph itself \cite{MARUSIC198169,TURNER1967136,MarSca}, its girth \cite{Eiben+Jaycay+Sparl+2019,Potocnik+Vidali+2022} etc...
    
  However, in the study of vertex-transitive graphs, it is often the case that vertex-transitivity itself (that is, the fact that any vertex can be mapped to any other via an automorphism) is not used in its full strength; instead, what is used is another weaker form of regularity regarding the distribution of cycles of a prescribed length throughout the graph, and that follows from vertex-transitivity. This observation motivated Poto\v{c}nik and Vidali to introduce the notion of girth-regularity (which generalises the concept of edge-girth-regularity of \cite{Jajcay+Kiss+Miklavic+2018}). Given a cubic graph, each vertex can be assigned a {\em signature} $(a,b,c)$, where the values $a$, $b$ and $c$ indicate the number of girth cycles containing each of the edges incident to the vertex. A graph is \emph{girth-regular} with signature $(a,b,c)$ if every vertex has signature $(a,b,c)$, a property that is automatically satisfied by vertex-transitive graphs.

     The notion of signature is a useful tool in the study and classification of cubic vertex-transitive graphs. For example, in \cite{CVTlargeorder}, where cubic vertex-transitive graphs with automorphisms of large order are classified. In \cite{Potocnik+Vidali+2019} and \cite{Potocnik+Vidali+2022} together, cubic vertex-transitive graphs of girth $g \leq 6$ were classified according to their signatures, extending the results of \cite{Eiben+Jaycay+Sparl+2019}. In this paper we continue this line of research, and further extend the results of Poto\v{c}nik and Vidali, by classifying all cubic vertex-transitive graphs of girth $7$. We show that only five signatures can occur for a cubic vertex-transitive graph of girth $7$, and that such a graph is either one of four exceptional graphs of small order, or belongs to one of three infinite families of graphs. 

   Besides the obvious advantage of having more comprehensive classification results for graphs, a better understanding of cubic vertex-transitive graphs of girth seven is important because of their relationship with an notable family of regular maps on surfaces, called Hurwitz maps. In \cite{Hurwits} Hurwitz showed that the automorphism group $G$ of a compact Riemann surface $S$ of genus $x \geq 2$ has order at most $84(x - 1)$, and that this bound is attained if and only if $S$ supports an orientably-regular map of type $\{7,3\}$ or $\{3,7\}$ whose group of orientation-preserving automorphisms is precisely $G$. Such a map is called a {\em Hurwitz map} and its skeleton (or the that of its dual) is a cubic vertex-transitive graph of girth at most $7$.

   Before we state the main theorem of this paper, we formally introduce, and briefly discuss, some of the notions appearing in it. 
   
\medskip
\noindent
{\sc Graphs.}
Throughout this paper, all graphs are finite and simple, unless stated otherwise.
In what follows, let $\Gamma$ be a graph. Then $V(\Gamma)$, $E(\Gamma)$, $A(\Gamma)$, and $\Aut(\Gamma)$ denote its vertex set, edge set, arc set and the automorphism group, respectively. An automorphism group $G \leq \Aut(\Gamma)$ is said to be \emph{vertex-transitive}, \emph{edge-transitive} and \emph{arc-transitive}, if the action of $G$ on the set $V(\Gamma)$ or the induced actions on $E(\Gamma)$, $A(\Gamma)$ are transitive.
Let $n\geq 3$ and let $\gamma = (v_0, \ldots, v_n)$ be a $(n+1)$-tuple of vertices of $\Gamma$ such that any two consecutive vertices are adjacent in $\Gamma$; we say that $\gamma$ is a \emph{walk} of length $n$. If $v_0 = v_n$ we say that $\gamma$ is \emph{closed} and if in addition any two vertices in $\{v_0, \ldots, v_{n-1}\}$ are distinct, then $\gamma$ is a \emph{cycle}. We consider two cycles to be the same if their vertices induce the same subgraph of $\Gamma$.

The precise definition of girth-regularity is the following. Let $\Gamma$ be a simple $d$-valent graph of girth $g$. For an edge $e \in E(\Gamma)$, we define $\epsilon(e)$ as the number of girth cycles that contain the edge $e$. Similarly, for a group of automorphisms $G$, we say that the signature of an $G$-edge orbit $\epsilon(e^G) = \epsilon(e)$. For a vertex $v \in V(\Gamma)$, let $\{e_1, \ldots, e_d\}$ be the set of edges incident to $v$, ordered in a way such that $\epsilon(e_1)\leq \ldots \leq \epsilon(e_d)$.
We say that the $d$-tuple $(\epsilon(e_1), \ldots, \epsilon(e_d))$ is the \emph{signature} of $v$. If every vertex in $\Gamma$ has the same signature $(\epsilon(e_1), \ldots, \epsilon(e_d))$, then $\Gamma$ is \emph{girth-regular} with signature $(\epsilon(e_1), \ldots, \epsilon(e_d))$. Note that if $\epsilon(e_1) = \epsilon(e_d)$, then $\Gamma$ is called \emph{edge-girth-regular} with parameters $(|V(\Gamma)|,d, g, \epsilon(e_1))$, where $g$ is the girth of $\Gamma$, see \cite{Jajcay+Kiss+Miklavic+2018}.

\medskip
\noindent
{\sc Maps.}
A map $M$ is an embedding of a graph $\Gamma$ on a compact surface $S$ with no boundary such that each connected component of $S-\Gamma$ is homeomorphic to an open disk.
There exist several equivalent definitions, but when interested only in embeddings of simple cubic graphs, many of the notions of the theory of maps can be defined in purely graph-theoretical terms. In this case, we  define a (trivalent) \emph{map} as a pair $M:= (\Gamma, \mathcal{C})$ where $\Gamma$ is a simple cubic graph and $\mathcal{C}$ is a set of cycles of $\Gamma$ such that every edge of $\Gamma$ is contained in precisely two cycles of $\mathcal{C}$. 
The vertex-, edge- and face-set of the $M$ are $V(M) = V(\Gamma)$, $E(M) = E(\Gamma)$ and $F(M) = \mathcal{C}$, respectively. If every walk in $\mathcal{C}$ has length $k$, we say that the map is of the type $\{k, 3\}$ ($3$ here is because the graph is cubic).
 A {\emph{flag}} of the map $M$ is a triple $(v,e,f)$, where $e \in E(M), v \in V(M), f \in F(M)$ and $v$ is incident to $e$, and $f$ contains $e$. The automorphism group of $M$ is denoted $\Aut(M)$ and consists of all the automorphisms of $\Gamma$ that preserve the set $\mathcal{C}$. We say that $M$ is a regular map, if $\Aut(M)$ acts transitively on the set of flags of $M$. If this action is transitive, it is also regular, hence the use of the term 'regular'. Similarly, if $\Aut(M)$ has the property that the
stabiliser of each vertex contains a cyclic group acting transitively on the neighbouring edges, and the stabiliser of each face contains a cyclic group acting transitively on the vertices of the face, then we say the map is \emph{rotary} (the term orientably-regular is also used to denote the same notion). Informally, a rotary map is one that possesses the maximum degree of rotational symmetry. In this case $\Aut(M)$ acts transitively on the vertices, arcs, faces, and incident pairs $(v,f)$  where $v$ is a vertex and $f$ is a face of $M$ containing $v$, but not necessarily on the flags (it may have up to $2$ orbits).

\medskip
\noindent
{\sc Dihedral schemes and truncations.} Let us now define truncations (with respect to dihedral schemes), which are a common concept in graph theory. For our definition we will follow \cite{Potocnik+Vidali+2019}, other and sometimes more general definitions appear in for example \cite{Eiben+Jaycay+Sparl+2019, Exoo+Jaycay+2012, Sachs+1963}, as well as some properties of truncations of graphs. Let $\Gamma$ be a graph, let $A(\Gamma)$ be the set of its arcs, and let $\leftrightarrow$ be an ireflexive symmetric relation on the set $A(\Gamma)$. 
Consider the graph $(A(\Gamma), \leftrightarrow)$ that has the vertex set $A(\Gamma)$ and with edge set $\{uv \mid u \leftrightarrow v\}$. 
Then $\leftrightarrow $ is a \emph{dihedral scheme} on $\Gamma$ if the connected components of the graph $(A(\Gamma), \leftrightarrow)$ are cycles such that each cycle contains precisely all the arcs with the same beginning. 
A \emph{truncation} of $\Gamma$ with respect to $\leftrightarrow$ is a simple graph that has the vertex set $A(\Gamma)$ and where two arcs $u, v  \in A(\Gamma)$ are adjacent when $u \leftrightarrow v$, or when $u$ and $v$ have the same underlying edge. Finally, let $\Aut(\Gamma, \leftrightarrow)$ be the group of all automorphisms of $\Gamma$ that preserve $\leftrightarrow$. We say that the dihedral scheme $\leftrightarrow$ is \emph{arc-transitive}, if $\Aut(\Gamma, \leftrightarrow)$ acts transitively on $A(\Gamma)$.

\medskip
We can now state the main theorem of this paper, where we let $\Pet(n,k)$ denote the generalized Petersen graph of order $2n$ and jump parameter $k$ (see \cite{Frucht+Graver+Watkins+1971} for a precise definition), and the graphs $\A(n)$ are Cayley graphs defined in \cref{sec:455-466-556}. 

\begin{theorem}\label{main_theorem}
    Let $\Gamma$ be a simple connected cubic vertex-transitive graph of girth $7$, then it is one of the following:
    \begin{enumerate}
        \item $\Gamma$ has signature $(0, 1, 1)$ and is isomorphic to the truncation of a $7$-regular graph $\Lambda$ (possibly with parallel edges) with respect to an arc transitive dihedral scheme $\leftrightarrow$.
        \item $\Gamma$ has signature $(2,2,2)$ and is the skeleton of a rotary map of type $\{7, 3\}$ embedded on a surface with the Euler characteristic $\chi = n( \frac{3}{7} - \frac{1}{2} )$.
        \item $\Gamma$ has signature $(4,4,4)$ and is the Coxeter graph.
        \item $\Gamma$ has signature $(4,4,6)$ and is isomorphic to the Cayley graph $\A(3n)$ for some $n \geq 3$.
        \item $\Gamma$ has signature $(4,5,5)$ and is one of the graphs $\Pet(13,5)$, $\Pet(15,4)$ and $\Pet(17,4)$.
    \end{enumerate}
\end{theorem}

It is worth mentioning that \cref{main_theorem} is 'almost' a characterization theorem. We will show that the graphs in items $(3)$, $(4)$ and $(5)$ are always connected cubic vertex-transitive graphs of girth $7$.  However, the graphs in items $(1)$ and $(2)$ are connected cubic vertex-transitive graphs but their girth might be strictly smaller (but never larger) than $7$. Nevertheless, cubic vertex-transitive graphs of girth $g\leq6$ are characterised in \cite{Potocnik+Vidali+2019, Potocnik+Vidali+2022}. Therefore, a graph in items $(1)$-$(5)$ is either a cubic vertex-transitive graph of girth $7$, or one of smaller girth appearing in \cite{Potocnik+Vidali+2019, Potocnik+Vidali+2022}.

To prove Theorem~\ref{main_theorem}, we will first determine a small set of possible signatures for a cubic vertex-transitive graph of girth $7$. In \cref{sec:possign}, we reduce the set of all possible signatures to a subset of seven elements. Graphs with signature $(0,1,1)$ or $(2,2,2)$ were first studied in \cite{Potocnik+Vidali+2019} and are discussed in \cref{sec:222}. Signature $(4,4,6)$ yields an infinite family of graphs and is analysed in \cref{sec:446}. In \cref{sec:455-466-556}, we work out the four remaining possible signatures, which yield only finitely many (or no) graphs.

\section{Possible signatures}
\label{sec:possign}

We start this section by stating the conditions for the signature of cubic girth-regular graphs that were proved in \cite{Potocnik+Vidali+2019}. 

\begin{lemma}[{\cite[Lemma 3.1, Lemma 3.2, Lemma 3.4]{Potocnik+Vidali+2019}}]\label{P1}
    If $(a, b, c)$ is the signature of a cubic girth-regular graph $\Gamma$ of girth $7$, then:
    \begin{enumerate}
        \item $a + b + c$ is even,
        \item If $a = 0$, then $(a, b, c) = (0,1,1)$,
        \item If $a \geq 1$, then $a + b > c$,
        \item $a + 4 \geq c$, and
        \item $a + 8 \geq b + c$.
    \end{enumerate}
\end{lemma}

We state another result that we will use to determine the possible signatures of cubic vertex-transitive graphs of girth $7$.
Nedela and \v Skoviera prove in \cite{Nedela+Skoviera+1995} that in a connected cubic graph, we cannot separate a cycle from another cycle by removing less than $g(\Gamma)$ edges from the graph.
More precisely, let $\Gamma$ be a connected graph. A subset $L \subseteq E(\Gamma)$ is \emph{cycle separating}, if $\Gamma - L$ is disconnected, and at least two of its components contain cycles. 
The \emph{cyclic connectivity} $\zeta(\Gamma)$ is the largest integer $k$ such that every cycle separating set has at least $k$ elements.

\begin{theorem}[{\cite[Theorem 17]{Nedela+Skoviera+1995}}]\label{P4}
    If $\Gamma$ is a connected cubic vertex-transitive graph, then $\zeta(\Gamma) = \text{girth}(\Gamma)$.
\end{theorem}

As a consequence of this theorem, we can see that the sum of parameters $a, b$ and $c$ is bounded by $18$:

\begin{proposition}\label{prop:P4}
    Let $\Gamma$ be a cubic vertex-transitive graph of girth $7$ and signature $(a,b,c)$. Then $a+b+c < 18$. 
\end{proposition}

\begin{proof}
    First assume that $|V(\Gamma)| \geq 28$.
    Let $v$ be a vertex in $\Gamma$, let $U \subseteq V(\Gamma)$ be the subset of all vertices at distance at most $3$ from $v$ and let $W = V(\Gamma) \setminus U$.
    Since $\Gamma$ has girth $7$, no two vertices at distance less than $3$ from $v$ can be neighbors or at distance $2$. Therefore $|U| = 22$ and $|W| \geq 6$. Let $\delta(U)$ be the set of edges that are incident to precisely one vertex in $U$. Suppose that $|\delta(U)| \leq 6$. Note that any girth-cycle through $v$ must contain only vertices in $U$, and thus the graph $\Gamma[U]$ contains a cycle. The graph $\Gamma[W]$ contains at least $6$ vertices and the average degree of a vertex in $\Gamma[W]$ is $\frac{1}{|W|}(\sum_{w \in W} \text{deg}(w)) = \frac{3 |W| - |\delta(U)|}{|W|} \geq 2$. Since $\Gamma[W]$ is a graph with average degree at least $2$, it contains a cycle.
    But then $\delta(U)$ is an edge-cut of $\Gamma$ of size less than $7$, such that that at least two of the parts of the disconnected graph $\Gamma - \delta(U)$ contain a cycle, which contradicts \cref{P4}.
    Hence $|\delta(U)| \geq 7$.
    
    Let $U' \subset U$ be the set of vertices at distance exactly $3$ from $v$. Note that each vertex in $U'$ is adjacent to exactly one vertex in $U\setminus U'$ (for otherwise there would a cycle of length shorter than $7$), and thus the number of edges joining a vertex in $U'$ with one in $U \setminus U'$ is precisely $12$. Moreover, the number of edges joining two vertices in $U'$ is $(a + b + c)/2$ since every girth cycle through $v$ must contain exactly one such edge. Thus
    \[ 36 = \sum_{w \in U'} \text{deg}(w) = |\delta(U)| + (a + b + c) + 12 . \]
    
It follows that $(a + b + c) \leq 17$. A computer search using the GAP package GraphSym \cite{GraphSym} shows that the result also holds if $\Gamma$ has less than $28$ vertices. 
See also \cite[Table 1]{Potocnik+Vidali+2022}. 
\end{proof}

\subsection{Graphs with multiple edge orbits} 

Most of the possible signatures do not contain only the same integer, which imposes multiple orbits on the set of edges. We examine such signatures. In the following $\lmulti \cdot \rmulti$ denotes a multiset.

\begin{de}
    Let $\Gamma$ be a graph, let $\O \subseteq \E(\Gamma)$ and let $c$ be a cycle in $\Gamma$. Let $\r(c, \O)$ be the number of edges in $c$ that are elements of $\O$. Similarly, for a vertex $v \in \V(\Gamma)$, let
    \[ \r(v, \O) = \lmulti  r(c, \O) \mid c \text{ is a girth-cycle that contains } v \rmulti  . \]
    If every vertex in $\Gamma$ has the same $\r$-number, we say that \[\r(\O) = \r(v, \O)\] is the $\r$-number of the set $\O$.
\end{de}

This notion is particularly useful when symmetry is involved. If $\Gamma$ is a $G$-vertex-transitive graph and $\O$ is a $G$-orbit of edges, then all vertices have the same $\r$-number. Note also that if the signature of $\Gamma$ is $(a_1, \ldots, a_t)$, then every vertex is contained in precisely in $\ell = (a_1, \ldots, a_t)/2$ girth-cycles, and thus the multiset $\r(\O)$ has size $\ell$. We will find some bounds and conditions for $\r(\O)$ that will help us understand the structure of $\Gamma$.

\begin{lemma}\label{lemma:sum}
    Let $\Gamma$ be a $G$-vertex-transitive graph of girth $g$ and let $\O$ be a $G$-edge-orbit. Then
    \[ \sum_{r_i \in \r(\O)} r_i = \frac{|\O|\epsilon(\O)g}{|\V(\Gamma)|}.\]
\end{lemma}

\begin{proof}
     First, let $\r(\O) = \lmulti  r_1, \ldots, r_k\rmulti  $ and let
    \begin{align*}
    \mathcal{A} = \{(c,e) \mid  e \in \O \text{ and $c$ is a girth-cycle containing $e$}\}.
    \end{align*}
    We will double count the size of the set $\mathcal{A}$.
     For $v \in \V(\Gamma)$, let
    \begin{align*}
    \mathcal{A}_v = \{(c,e) \in \mathcal{A}\mid \text{the vertex $v$ is contained in $c$} \}.
    \end{align*}
    Observe that $|\mathcal{A}_v|= r_1 + \ldots + r_k$. The pair $(c, e) \in \mathcal{A}$ appears in $\mathcal{A}_v$ precisely if $v$ is a vertex of $c$. Hence
    \begin{align*}
        \sum_{u \in V(\Gamma)} |\mathcal{A}_u| = g |\mathcal{A}|.
    \end{align*}
    Since $\Gamma$ is a vertex transitive graph, it follows that $|\mathcal{A}_v| = |\mathcal{A}_u|$ for all $u \in \V(\Gamma)$ and thus
    \begin{align*}
        |\mathcal{A}| = (r_1 + \ldots + r_k)\frac{|V(\Gamma)|}{g}.
    \end{align*}
    Another way to count the pairs $(c,e)$ is to observe that every edge $e \in \O$ is contained in precisely $\epsilon(\O)$ girth cycles. Hence 
    \begin{align*}
        |\mathcal{A}| = \sum_{e \in \O} \epsilon(e) =|\O|\epsilon(\O),
    \end{align*}
    and thus 
    \[
        r_1 + \ldots + r_k = \frac{|\O|\epsilon(\O) g}{|\V(\Gamma)|}.
    \qedhere \]
 \end{proof}
\medskip

\begin{lemma}\label{lem: repeatedodd}
    Let $\Gamma$ be a  cubic vertex-transitive girth-regular graph of odd girth $g$ with signature $(a,b,c)$ and let $s \in \{a,b,c\}$. If $s$ is an odd integer, then it appears exactly twice in the signature. In that case, any vertex-transitive group of automorphisms $G$ of $\Gamma$ has precisely two orbits for its action on the edges of $\Gamma$.
\end{lemma}

\begin{proof}
    Suppose that $s \in \{a, b, c\}$ is an odd integer and let $e$ be an edge in $\Gamma$ such that $\epsilon(e) = s$. Suppose that there exists a vertex transitive automorphism group $G$ of $\Gamma$ such that $|e^G| = \frac{|E(\Gamma)|}{3}$, and let $v$ be a vertex that is incident to $e$. Then by Lemma~\ref{lemma:sum}, it follows that 
    \begin{align}\label{eq1}
        r_1 + \ldots + r_k = \frac{gs|E(\Gamma)|}{3 |V(\Gamma)|} = \frac{g}{2}s,
    \end{align}
    where $\lmulti r_1, \ldots, r_k \rmulti = \r(e^G)$ are as in the statement of Lemma~\ref{lemma:sum}. Since the left hand side of Equation~\ref{eq1} is an integer, and both $s$ and $g$ are odd integers, we arrive to a contradiction.

    It follows that $e^G$ contains at least $\frac{2}{3}|\E(\Gamma)|$ edges, and $G$ has at most $2$ orbits on edges. Since not all three parameters $a, b$ and $c$ can be odd by Lemma~\ref{P1}, $G$ can not be transitive on the set of edges. Therefore it must have exactly two orbits for its action on $E(\Gamma)$ and $s$ appears exactly twice in the signature.
\end{proof}

By \cref{lem: repeatedodd}, \cref{P1} and \cref{prop:P4}, we see that the possible signatures of a cubic vertex-transitive graph of girth $7$ are: $(0, 1, 1), (2, 2, 2), (2, 3, 3), (2, 4, 4), (2, 5, 5), (3, 3, 4), (4, 4, 4), (4, 4, 6), (4, 5, 5), (4, 6, 6) \\ \text{ and }  (5, 5, 6)$. We will eliminate some additional signatures with the following lemma and its corollaries, before we go into more in-depth analysis in the following sections.

\begin{lemma} \label{lemma:multiplicity}
   Let $\Gamma$ be a cubic $G$-vertex-transitive graph of girth $7$ with signature $(a,b,c)$, and let $\O$ be a $G$-edge orbit such that $|\O| = |\E(\Gamma)|/3$. Then if either of the integers $1, 2$ or $3$ is contained in the multiset $\r(\O)$, it is contained in it at least four times. If $0$ is contained in $\r(\O)$, it is contained exactly once. 
\end{lemma}

\begin{proof}
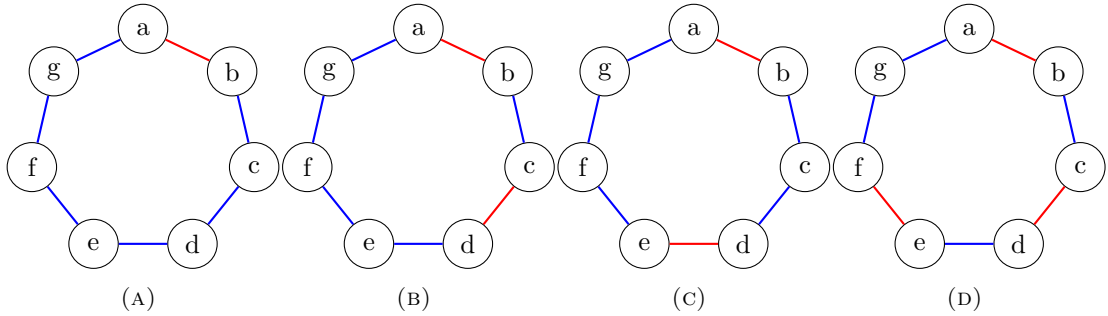
\begin{figure}[h]
  \begin{subfigure}[b]{0.22\textwidth} 
    \begin{tikzpicture}
\graph [clockwise=7, nodes={draw, circle, minimum size=0.65cm}, simple, radius=1.5 cm] { a, b, c, d, e, f, g;
a -- [red, thick] b;
b -- [blue, thick] c;
c -- [blue,thick] d;
d -- [blue,thick] e;
e -- [blue,thick] f;
f -- [blue,thick] g;
g -- [blue,thick] a;
};
\end{tikzpicture}
\caption[]{}
\label{figure:c1}
\end{subfigure}
  \begin{subfigure}[b]{0.22\textwidth} 
    \begin{tikzpicture}
\graph [clockwise=7, nodes={draw, circle, minimum size=0.65cm}, simple, radius=1.5 cm] { a, b, c, d, e, f, g;
a -- [red, thick] b;
b -- [blue,thick] c;
c -- [red, thick] d;
d -- [blue,thick] e;
e -- [blue,thick] f;
f -- [blue,thick] g;
g -- [blue,thick] a;
};
\end{tikzpicture}
\caption[]{}
\label{figure:c2}
\end{subfigure}
  \begin{subfigure}[b]{0.22\textwidth} 
    \begin{tikzpicture}
\graph [clockwise=7, nodes={draw, circle, minimum size=0.65cm}, simple, radius=1.5 cm] { a, b, c, d, e, f, g;
a -- [red, thick] b;
b -- [blue,thick] c;
c -- [blue,thick] d;
d -- [red, thick] e;
e -- [blue,thick] f;
f -- [blue,thick] g;
g -- [blue,thick] a;
};
\end{tikzpicture}
\caption[]{}
\label{figure:c3}
\end{subfigure}
  \begin{subfigure}[b]{0.22\textwidth} 
    \begin{tikzpicture}
\graph [clockwise=7, nodes={draw, circle, minimum size=0.65cm}, simple, radius=1.5 cm] { a, b, c, d, e, f, g;
a -- [red, thick] b;
b -- [blue,thick] c;
c -- [red, thick] d;
d -- [blue,thick] e;
e -- [red, thick] f;
f -- [blue,thick] g;
g -- [blue,thick] a;
};
\end{tikzpicture}
\caption[]{}
\label{figure:c4}
\end{subfigure}
\caption[]{Possible $7$-cycles with one, two or three edges from $\O$.}
\label{fig:cycles}
\end{figure}
We will call the edges form $\O$ red edges, and the other edges blue.
Suppose first that $1 \in \r(\O)$. Then there exists a $7$-cycle $\gamma = (a,b,c,d,e,f,g,a)$ in $\Gamma$ that contains $v$ and exactly one red edge, we can assume this is the edge $ab$ (see \cref{figure:c1}). Take the set of vertices $\{a, c, d, e\}$ and note that the vertex $a$ is incident to one red edge, vertices $c$, $d$ and $e$ are at distance $1$, $2$ and $3$ from a red edge, respectively. Since $G$ preserves the set of red edges, there exists no automorphism preserving $\gamma$ that maps one of $\{a, c, d, e\}$ into another one of the vertices $\{a, c, d, e\}$.
Hence the automorphisms that map vertices $a, c, d$ and $e$ to $v$, map the cycle $\gamma = (a,b,c,d,e,f,g,a)$ to four distinct $7$-cycles that contain $v$ and exactly one red edge. Hence the multiplicity of $1$ in $\r(\O)$ is at least $4$.

Similarly, if $2 \in \r(\O)$, then there exists a $7$-cycle in $\Gamma$ that contains $v$ and two edges from the orbit $\O$ as in one of the Figures \ref{figure:c2} or \ref{figure:c3}, where edges from $\O$ are colored red. In the first case, take the set $\{a,b, e,f\}$, and in the second take the set $\{a, b, c, f\}$. Similarly as before, it follows that the multiplicity of $2$ in $\r(\O)$ is at least $4$.

If $3 \in \r(\O)$, then there exists a cycle as in \cref{figure:c4} that contains $v$ and three red edges. By considering the set $\{a, b, c, g\}$, we can show that the multiplicity of $3$ in $\r(\O)$ is at least $4$.

Since $\Gamma - \O$ is a $2$-valent graph, it is a union of blue cycles, hence there exists at most one such $7$-cycle containing a given vertex if $\Gamma$. It follows that $0$ appears in $\r(\O)$ at most once.
\end{proof}

From the proof of \cref{lemma:multiplicity} we can also see the following observation, that we will need in \cref{sec:455-466-556}.

\begin{remark}\label{rem:2red}
    Let $\Gamma$ be a cubic vertex-transitive graph of girth $7$, and let $\O$ be an edge orbit such that every vertex is incident to one edge from $\O$. Then if a vertex $v$ is contained in one cycle as in \cref{figure:c2}, it is contained in at least $4$ such cycles (the edges of the cycle are red-blue-red-blue-blue-blue-blue). If it is contained in one cycle as in \cref{figure:c3}, it is contained in at least $4$ cycles as \cref{figure:c3} (the edges of the cycle are red-blue-blue-red-blue-blue-blue).
\end{remark}

\cref{lemma:sum} and \cref{lemma:multiplicity} (or rather, \cref{cor:condition}, which summarizes them) will be used multiple times, thus we state, for convenience, the following condition for a pair $(\Gamma, \O)$, where $\Gamma$ is a cubic vertex-transitive graph and $\O$ is an edge-orbit under $\Aut(\Gamma)$.

\begin{condition}\label{CONDITION}
 There exist integers $r_1\leq \ldots\leq r_{\frac{(a + b+c)}{2}} \in \{0,1,2,3\}$ such that 
 \begin{enumerate}
     \item $r_1 + \dots + r_{\frac{(a + b+c)}{2}} = 7 \epsilon(\O)/2$;
     \item $r_2 > 0$;
     \item For every $r_i \neq 0$ there exist at least three indices $j,k,\l$ such that $i,j,k,\l$ are pairwise distinct, and $r_i = r_j = r_k = r_\ell$.
  \end{enumerate}
 \end{condition}

 \begin{cor}
 \label{cor:condition}
     Let $\Gamma$ be a cubic vertex-transitive graph of girth $7$. If there exist an edge orbit $\O$ such that $|\O| = |\E(\Gamma)/3|$, then the integers from $\r(\O)$: $r_1\leq \ldots\leq r_{\frac{(a + b+c)}{2}}$ satisfy \cref{CONDITION} for $(\Gamma, \O)$.
  \end{cor}

\begin{cor}
    Let $\Gamma$ be a   cubic vertex-transitive girth-regular graph of girth $7$ with signature $(a,b,c)$. Then $(a, b, c) \notin \{(2,3,3), (2,4,4), (2,5,5), (3,3,4)\}$.
\end{cor}

\begin{proof}
    Suppose for a contradiction that $\Gamma$ is a vertex-transitive cubic graph of girth $7$ and signature $(a,b,c) \in \{(2,3,3), (2,4,4), (2,5,5), (3,3,4)\}$. In each of the signatures, one of the integers is different from the others, hence the edges with this signature form an edge orbit $\O$. We can check that $\Gamma$ and $\O$ do not satisfy \cref{CONDITION}, which is a contradiction by \cref{cor:condition}.
\end{proof}

We are left with the signatures $(0, 1, 1), (2, 2, 2), (4, 4, 4), (4, 4, 6), (4, 5, 5), (4, 6, 6)$ and $(5, 5, 6)$. They will be examined in the next sections.

\section{Signatures $(0,1,1)$ and $(2,2,2)$}
\label{sec:222}

Graphs with signature $(0,1,1)$ and $(2,2,2)$ were studied in \cite{Potocnik+Vidali+2019} for arbitrary girth. In particular, those with signature $(0,1,1)$ can always be obtained from arc-transitive dihedral schemes on $7$-valent graphs. Every such dihedral scheme where the graph has girth at least $4$ (as well as lower girth if the dihedral scheme is appropriate) will give us a cubic vertex-transitive graph of girth $7$. Since not much is known about $7$-valent arc-transitive graphs, this family is large and difficult to understand.

\begin{theorem}[{\cite[Theorem 3.6]{Potocnik+Vidali+2019}}]\label{P2}
    If $\Gamma$ is a cubic girth-regular graph of girth $g$ with signature $(0,1,1)$, then $\Gamma$ is isomorphic to the truncation of a $g$-valent graph $\Lambda$ (possibly with parallel edges) with respect to a dihedral scheme $\leftrightarrow$. Moreover, if $\Gamma$ is vertex-transitive, then the dihedral scheme $\leftrightarrow$ is arc-transitive.
\end{theorem}

On the other hand, cubic graphs of girth $7$ and signature $(2,2,2)$ are skeletons of maps of type $\{7,3\}$.

\begin{theorem}[{\cite[Theorem 3.11]{Potocnik+Vidali+2019}}]\label{P3}
    Let $\Gamma$ be a connected cubic girth-regular graph of girth $g$ with $n$ vertices and signature $(2,2,2)$. Then $g$ divides $3n$ and $\Gamma$ is the skeleton of a map of type $\{g, 3\}$ embedded on a surface with the Euler characteristic $$\chi = n\biggl ( \frac{3}{g} - \frac{1}{2} \biggr ).$$
    Moreover, every automorphism of $\Gamma$ extends to an automorphism of the map. In particular,
if $\Gamma$ is vertex-transitive, so is the map.
\end{theorem}

We can say even more about symmetry properties of the map in \cref{P3} when $\Gamma$ is vertex-transitive and has girth $7$.

\begin{lemma}\label{lem:rotary7}
    Let $\Gamma$ be a connected cubic vertex-transitive graph of girth $7$ and signature $(2,2,2)$. Then $\Gamma$ is the skeleton of a rotary map of type $\{7,3\}$.
\end{lemma}

\begin{proof}
Let $\C$ be the set of girth cycles of $\Gamma$ and observe that $M := (\Gamma,\C)$ is a map of type $\{7,3\}$. We will show that $M$ is rotary. That is, that for every vertex, there is an automorphism fixing it while cyclically permuting its neighbors, and that for every face, there is an automorphism fixing it set-wise while cyclically permuting its vertices.

First, we will show that $G := \Aut(\Gamma)$ is edge-transitive. Suppose to the contrary, that $G$ has more than one edge-orbit. Then, necessarily one of these orbits has precisely $|\E(\Gamma)|/3$ edges. However, such an orbit (together with $\Gamma$) cannot satisfy \cref{CONDITION}, contradicting \cref{cor:condition}. Therefore $G$ must be edge-transitive, and since $G$ is vertex-transitive and $\Gamma$ has odd valency, $G$ must be arc-transitive as well. Then for every $v \in \V(\Gamma)$, the stabilizer $G_v$ acts on the neighbors of $v$ as $C_3$ or $\Dih(3)$. It follows that for every vertex, there is an automorphism fixing it and cyclically permuting its neighbors.

Now, to show that the condition on the faces hold, we will focus on the $2$-arcs of $\Gamma$. 
In particular, we will show that given two distinct $2$
-arcs $a$ and $a'$, there exists an automorphism mapping $a$ to either $a'$ or to its inverse $(a')^{-1}$. That is clearly satisfied if $G$ is $2$-arc-transitive, so we may assume that $G$ has two orbits on $2$-arcs, which we denote $\O_1$ and $\O_2$. since every $1$-arc is the beginning of two distinct $2$-arcs, we see that $\O_1$ and $\O_2$ are of equal size. Let $(u, v, w) \in \mathcal{O}_1$ and suppose that its inverse, $(w,v,u),$ is also in $\mathcal{O}_1$. Then for all $2$-arcs in $\mathcal{O}_1$, their inverse arc is also an element of $\mathcal{O}_1$. Since $\O_2$ contains all the arcs not in $\O_1$, the same is true for $\O_2$. Let $z \notin \{u,w\}$ be the third neighbor of $v$. Then, there are six $2$-arcs having $v$ as its middle vertex: $(u,v,w)$, $(u,v,z)$, $(z,v,w)$ and their inverses. Without loss of generality we can assume that $(u,v,z) \in \mathcal{O}_1$. But since every $2$-arc and its inverse are in the same orbit, we have that four of the $2$-arcs centered at $v$ belong to $\O_1$ while only two belong to $O_2$. By the vertex-transitivity of $G$ we have that for each vertex, four of the $2$-arcs centered at it belong to $\O_1$, implying that $\O_1$ is twice the size of $\O_2$, contradicting that both orbits have the same size. The contradiction stems from the assumption that there is an automorphism mapping $(u,v,w)$ to its inverse. It follows that a $2$-arc and its inverse belong to distinct orbits. This, in turn, means that given two distinct $2$-arcs $a$ and $a'$, there exists an automorphism mapping $a$ to either $a'$ or to its inverse.

Now, observe that since $\Gamma$ has signature $(2,2,2)$, every $2$-arc belongs to precisely one cycle in $\C$ (or one face of $M$, if you will). This means that each $2$-arc uniquely defines a face, and that mutually inverse $2$-arcs define the same face. Thus, $G$ must act transitively on $\C$, the set of faces of $M$, and given a face $F \in \C$ the stabilizer $G_F$ acts transitively on the vertices of $F$. Since each face is a $7$-cycle, $G_F$ must act on $F$ like a transitive subgroup of $\Dih(7)$. Since the only two transitive subgroups of $\Dih(7)$ are either $C_7$ or $\Dih(7)$ itself, we see that $G_F$ admits an automorphism cyclically permuting the vertices of $F$. 

We conclude that $M$ is a rotary map.
\end{proof}

\section{Signature $(4,4,6)$}
\label{sec:446}
\begin{figure}[h]
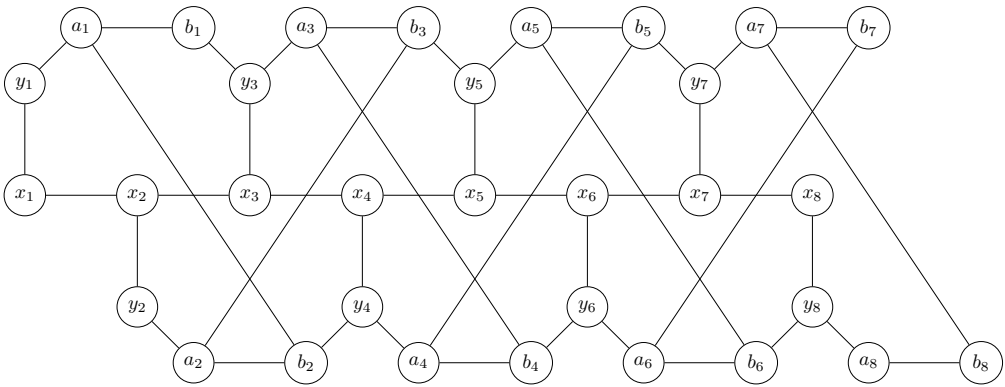

  \centering
\begin{resizedtikzpicture}{\textwidth}
	\begin{pgfonlayer}{nodelayer}
		\node [style=krog] (3) at (-4.5, 0.5) {$x_1$};
		\node [style=krog] (4) at (-2.5, 0.5) {$x_2$};
		\node [style=krog] (11) at (-4.5, 2.5) {$y_1$};
		\node [style=krog] (12) at (-2.5, -1.5) {$y_2$};
		\node [style=krog] (19) at (-3.5, 3.5) {$a_1$};
		\node [style=krog] (20) at (-1.5, 3.5) {$b_1$};
		\node [style=krog] (39) at (-0.5, 0.5) {$x_3$};
		\node [style=krog] (40) at (1.5, 0.5) {$x_4$};
		\node [style=krog] (42) at (-0.5, 2.5) {$y_3$};
		\node [style=krog] (43) at (1.5, -1.5) {$y_4$};
		\node [style=krog] (45) at (0.5, 3.5) {$a_3$};
		\node [style=krog] (46) at (2.5, 3.5) {$b_3$};
		\node [style=krog] (47) at (-1.5, -2.5) {$a_2$};
		\node [style=krog] (48) at (0.5, -2.5) {$b_2$};
		\node [style=krog] (51) at (3.5, 0.5) {$x_5$};
		\node [style=krog] (52) at (5.5, 0.5) {$x_6$};
		\node [style=krog] (54) at (3.5, 2.5) {$y_5$};
		\node [style=krog] (55) at (5.5, -1.5) {$y_6$};
		\node [style=krog] (57) at (4.5, 3.5) {$a_5$};
		\node [style=krog] (58) at (6.5, 3.5) {$b_5$};
		\node [style=krog] (59) at (2.5, -2.5) {$a_4$};
		\node [style=krog] (60) at (4.5, -2.5) {$b_4$};
		\node [style=krog] (63) at (7.5, 0.5) {$x_7$};
		\node [style=krog] (64) at (9.5, 0.5) {$x_8$};
		\node [style=krog] (66) at (7.5, 2.5) {$y_7$};
		\node [style=krog] (67) at (9.5, -1.5) {$y_8$};
		\node [style=krog] (69) at (8.5, 3.5) {$a_7$};
		\node [style=krog] (70) at (10.5, 3.5) {$b_7$};
		\node [style=krog] (71) at (6.5, -2.5) {$a_6$};
		\node [style=krog] (72) at (8.5, -2.5) {$b_6$};
		\node [style=krog] (73) at (10.5, -2.5) {$a_8$};
		\node [style=krog] (74) at (12.5, -2.5) {$b_8$};
		
		\node [style=none] (102) at (-5, 3) {};
		\node [style=none] (103) at (-5.25, 0.5) {};
		\node [style=none] (104) at (-4.5, -1) {};
		\node [style=none] (105) at (-3, -2) {};
		
		\node [style=pika] (120) at (13.75, 0.5) {};
		\node [style=pika] (121) at (14.25, 0.5) {};
		\node [style=pika] (122) at (14.75, 0.5) {};
            \node [style=pika] (120) at (-5.75, 0.5) {};
		\node [style=pika] (121) at (-6.25, 0.5) {};
		\node [style=pika] (122) at (-6.75, 0.5) {};
	\end{pgfonlayer}
	\begin{pgfonlayer}{edgelayer}
		\draw (11) to (19);
		\draw (19) to (20);
		\draw (20) to (42);
		\draw (39) to (4);
		\draw (4) to (3);
		\draw (11) to (3);
		\draw (12) to (4);
		\draw (12) to (47);
		\draw (19) to (48);
		\draw (47) to (48);
		\draw (42) to (45);
		\draw (45) to (46);
		\draw (46) to (54);
		\draw (51) to (40);
		\draw (40) to (39);
		\draw (42) to (39);
		\draw (48) to (43);
		\draw (43) to (40);
		\draw (43) to (59);
		\draw (47) to (46);
		\draw (45) to (60);
		\draw (59) to (60);
		\draw (54) to (57);
		\draw (57) to (58);
		\draw (58) to (66);
		\draw (63) to (52);
		\draw (52) to (51);
		\draw (54) to (51);
		\draw (60) to (55);
		\draw (55) to (52);
		\draw (55) to (71);
		\draw (59) to (58);
		\draw (57) to (72);
		\draw (71) to (72);
		\draw (66) to (69);
		\draw (69) to (70);
		\draw (64) to (63);
		\draw (66) to (63);
		\draw (72) to (67);
		\draw (67) to (64);
		\draw (67) to (73);
		\draw (71) to (70);
		\draw (69) to (74);
		\draw (73) to (74);
	\end{pgfonlayer}
\end{resizedtikzpicture}

\caption[The graph $\A(n)$]{A section of the graph $\A(n)$.}
\label{figure:An}
\end{figure}

In this section we deal with 
vertex-transitive graphs of girth $7$ and signature $(4,4,6)$. We show that they form an infinite family of Cayley graphs.
For $n \geq 8$, define $\A(n)$ as the graph with the vertex-set 
\begin{align*}
    \V(\A(n)) &= \{x_i \mid i \in \ZZ_n \} \cup \{y_i \mid i \in \ZZ_n \} \cup \{a_i \mid i \in \ZZ_n \} \cup \{b_i \mid i \in \ZZ_n \},
\end{align*}
and edge-set
\begin{align*}
    \E(\A(n)) &= \{x_i x_{i+1},\ x_i y_i,\ y_i a_i,\ a_i b_i,\ b_i y_{i+2},\  a_i b_{i+1} \mid i \in \ZZ_n  \}.
\end{align*}
A part of the graph $\A(n)$, with $n\geq 8$, is shown in \cref{figure:An}.

\begin{lemma}\label{lem:446_1}
If $n \geq 8$, then $\A(n)$ has girth $7$ and signature $(4,4,6)$. It is vertex-transitive if and only if $n$ is divisible by $3$.
\end{lemma}

\begin{proof}
    Let $G = \Aut(\A(n))$ and notice that the mapping $s_i \mapsto s_{i+1}$, for $s \in  \{a,b,x,y\}$, is an automorphism of $\A(n)$. It follows that the vertices of the form $x_i$ for $i \in \ZZ_n$ are contained in the same $G$-orbit. Similarly, each of the sets $\{y_i \mid i \in\ \ZZ_n\}$, $\{a_i \mid i \in\ \ZZ_n\}$ and $\{b_i \mid i \in\ \ZZ_n\}$ is contained in some $G$-orbit. By the same argument all edges of the form $x_i x_{i+1}$ belong to the same orbit. An analogous claim can be made for the edges of each of the following forms: $x_i y_i$, $y_i a_i$, $a_i b_i$, $b_i y_{i+2}$, and $a_i b_{i+1}$. By the definition of $A(n)$ (or by inspecting \cref{figure:An}), we see that $\epsilon(x_i x_{i+1}) = 4$, $\epsilon(x_i y_i) = 6$, $\epsilon(y_i a_i) =  4  $, $\epsilon(a_i b_i) = 6$, $\epsilon(b_i y_{i+2}) = 4$ and $\epsilon(a_i b_{i+1}) = 4$, and thus the signature of every vertex is $(4,4,6)$. 

    Now, suppose that $n$ is not divisible by $3$. Since the edges $x_iy_i$ and $a_ib_i$ have signature $6$, they form a $1$-factor $F$ in $\A(n)$ that is preserved by $G$. Then $\A(n) - F$ is a union of cycles. Clearly, $\gamma_1 = (x_1, x_2, \ldots, x_n, x_1)$ and  $\gamma_2 = (y_1, a_1, b_2, y_4, a_4, b_5, y_7, a_7, \ldots)$ are two such cycles. However, $\gamma_1$ has length $n$, and $\gamma_2$ has length $3n$ (since $n$ is not divisible by $3$). It follows that $\A(n)$ is not vertex transitive.
    
    For the converse, assume that $n = 3k$ for an integer $k \in \NN$. The map $\tau$ with 
    \begin{align*}
    \tau(x_i) = 
    \begin{cases}
    b_{i-1} & \text{if }i \equiv 0 \pmod 3\\
    y_i & \text{if }i \equiv 1 \pmod 3 \\
    a_{i-1} & \text{if }i \equiv 2 \pmod 3
    \end{cases}\ \ \ 
    \tau(y_i) = 
    \begin{cases} 
    a_{i-1} & \text{if }i \equiv 0 \pmod 3\\
    x_i & \text{if }i \equiv 1 \pmod 3 \\
    b_{i-1} & \text{if }i \equiv 2 \pmod 3
    \end{cases} \\
    \tau(a_i) = 
    \begin{cases}
    b_{i} & \text{if }i \equiv 0 \pmod 3\\
    x_{i+1} & \text{if }i \equiv 1 \pmod 3 \\
    y_{i+1} & \text{if }i \equiv 2 \pmod 3
    \end{cases}\ \ \ 
    \tau(b_i) = 
    \begin{cases}
    a_{i} & \text{if }i \equiv 0 \pmod 3\\
    y_{i+1} & \text{if }i \equiv 1 \pmod 3 \\
    x_{i+1} & \text{if }i \equiv 2 \pmod 3
    \end{cases}
    \end{align*}
    
    is an automorphism of $\A(n)$. Since each of the sets $\{x_i \mid i \in\ \ZZ_n\}$, $\{y_i \mid i \in\ \ZZ_n\}$, $\{a_i \mid i \in\ \ZZ_n\}$ and $\{b_i \mid i \in\ \ZZ_n\}$ is contained in some $G$-orbit, this is enough to conclude that  $\A(n)$ is vertex-transitive. 
\end{proof}

\begin{lemma}\label{lem:446_2}
    If $\Gamma$ is a connected vertex-transitive graph with girth $7$ and signature $(4,4,6)$, then $\Gamma$ is isomorphic to $\A(\frac{k}{4})$, where $k = |V(\Gamma)|$.
\end{lemma}

\begin{proof}

    Let $G = \Aut(\Gamma)$ and let $\O$ denote the orbit of edges with signature $6$. For convenience, we will call the edges in $\O$ \emph{red} edges, and those in $E(\Gamma) \setminus \O$ \emph{blue} edges. Then, every vertex $v \in V(\Gamma)$ is incident to a unique red edge, and two blue edges, and is contained in exactly $7$ girth cycles. Let $r_1, \ldots, r_7$ denote the number of red edges contained in each of those cycles. The only integers that satisfy \cref{CONDITION} and \cref{cor:condition} are
    $r_i = 3$ for $i = 1, \ldots, 7$, and thus each girth cycle in $\Gamma$ contains precisely $3$ red edges. 
    
    Since $G$ preserves the set $E(\Gamma) \setminus \O$ of blue edges, and since $\Gamma$ is vertex-transitive, it follows that for every two consecutive blue edges  $vu,uw$, there exists a girth cycle that contains the walk $(v,u,w)$. We will use this fact often.

    The graph $\Gamma - \O$, induced by the blue edges, is a disjoint union of cycles. Let $(x_1, x_2, \ldots, x_n, x_1)$ be one of those cycles. Each of the vertices $x_1, \ldots, x_n$ is incident to exactly one red edge. Suppose that one such edge is of the form $x_ix_j$ for some distinct $i,j \in \{1, \ldots, n\}$. Since $G$ preserves the set of blue edges, we see that $\Delta = \{x_1, \ldots, x_n\}$ is a block of imprimitivity for this action. Since those vertices form a cycle, $G_{\{\Delta\}}^\Delta$ is a permutation group isomorphic to a subgroup of the dihedral group $\text{D}_n$ of order $2n$. Since $G$ is vertex-transitive, this subgroup is transitive on $\{x_1, \ldots, x_n\}$. Every transitive subgroup of $\Dih(n)$ must contain a subgroup isomorphic to $2\ZZ_n$ generated by a two-step rotation, hence $x_{i + 2\l} x_{j + 2\l} \in \O$ for every integer $\l$. We see that $\Gamma$ contains a cycle $(x_{i-2}, x_{i-1}, x_i, x_j, x_{j-1}, x_{j-2}, x_{i-2})$ of length $6$, which is a contradiction.

    It follows that for every $i \in \ZZ_n$ there exist a vertex $y_i \notin \{x_1, \ldots, x_n\}$ such that $x_iy_i$ is a red edge. Since each vertex in $\Gamma$ is incident to precisely one red edge, we see that all the vertices $y_i$ are distinct. So far we have shown that $\Gamma$ contains a subgraph as in \cref{figure:An_partial}, where the edges from $\O$ are colored red.

\begin{figure}[h]
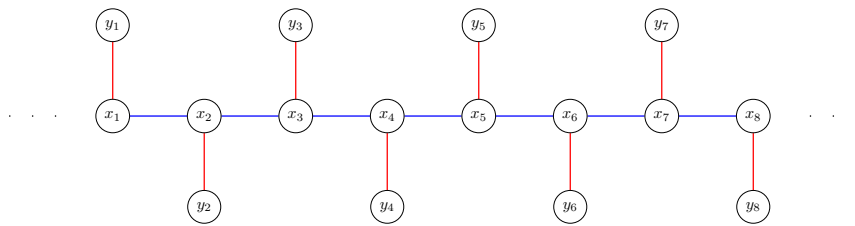

  \centering
\begin{resizedtikzpicture}{0.7\textwidth}
	\begin{pgfonlayer}{nodelayer}
		\node [style=krog] (3) at (-4.5, 0.5) {$x_1$};
		\node [style=krog] (4) at (-2.5, 0.5) {$x_2$};
		\node [style=krog] (11) at (-4.5, 2.5) {$y_1$};
		\node [style=krog] (12) at (-2.5, -1.5) {$y_2$};
		\node [style=krog] (39) at (-0.5, 0.5) {$x_3$};
		\node [style=krog] (40) at (1.5, 0.5) {$x_4$};
		\node [style=krog] (42) at (-0.5, 2.5) {$y_3$};
		\node [style=krog] (43) at (1.5, -1.5) {$y_4$};
		\node [style=krog] (51) at (3.5, 0.5) {$x_5$};
		\node [style=krog] (52) at (5.5, 0.5) {$x_6$};
		\node [style=krog] (54) at (3.5, 2.5) {$y_5$};
		\node [style=krog] (55) at (5.5, -1.5) {$y_6$};
		\node [style=krog] (63) at (7.5, 0.5) {$x_7$};
		\node [style=krog] (64) at (9.5, 0.5) {$x_8$};
		\node [style=krog] (66) at (7.5, 2.5) {$y_7$};
		\node [style=krog] (67) at (9.5, -1.5) {$y_8$};
		\node [style=pika] (120) at (10.75, 0.5) {};
		\node [style=pika] (121) at (11.25, 0.5) {};
		\node [style=pika] (122) at (11.75, 0.5) {};
            \node [style=pika] (120) at (-5.75, 0.5) {};
		\node [style=pika] (121) at (-6.25, 0.5) {};
		\node [style=pika] (122) at (-6.75, 0.5) {};
	\end{pgfonlayer}
	\begin{pgfonlayer}{edgelayer}
		\draw [blue, thick](39) to (4);
		\draw [blue, thick] (4) to (3);
		\draw[red, thick] (11) to (3);
		\draw[red, thick] (12) to (4);
		\draw [blue, thick] (51) to (40);
		\draw [blue, thick] (40) to (39);
		\draw[red, thick] (42) to (39);
		\draw[red, thick] (43) to (40);
		\draw [blue, thick] (63) to (52);
		\draw [blue, thick] (52) to (51);
		\draw[red, thick] (54) to (51);
		\draw[red, thick] (55) to (52);
		\draw [blue, thick] (64) to (63);
		\draw[red, thick] (66) to (63);
		\draw[red, thick] (67) to (64);
	\end{pgfonlayer}
\end{resizedtikzpicture}

\caption[\ ]{\ Vertices $x_1, \ldots, x_n$ and $y_1, \ldots, y_n$ in $\Gamma$. }
\label{figure:An_partial}
\end{figure}

Now since $x_i x_{i+1}$ and $x_{i+1} x_{i+2}$ are blue edges for each $i$, there exists a $7$-cycle that contains the walk $(x_i, x_{i+1}, x_{i+2})$. Since this $7$-cycle must contain three (non-consecutive) red edges, there exists a path $(y_i, a_i, b_i, y_{i+2})$ for each $i \in \ZZ_n$, and for some $a_i, b_i \in V(\Gamma)$ such that $a_i b_i$ is red, and both $y_ia_i$ and $b_i y_{i+2}$ are blue.

We claim that for every $i \in \ZZ_n$, the vertices $a_i$ and $b_i$ are distinct from vertices $\{x_1, \ldots, x_n, y_1, \ldots, y_n\}$, that $a_i \neq a_j$, $b_i \neq b_j$ for $i \neq j$ and that $b_i \neq a_j$ for any $i,j \in \ZZ_n$. First note that for every $i$, the vertices $a_i$ and $b_i$ are distinct from the vertices in $\{x_1, \ldots, x_n\}$, since no $x_i$ is adjacent to a $y_j$ through a blue edge. 
Similarly for all $i \in \ZZ_n$, the vertices $a_i$ and $b_i$ are distinct from the vertices in $\{y_1, \ldots, y_n\}$, since the edge $a_ib_i$ is red, but for every $k \in \{1, \ldots, n\}$, the vertex $y_k$ is already incident to a red edge $x_ky_k$, and $a_i, b_i \notin\{x_1, \ldots, x_n\}$.
Suppose now for a contradiction that there exist $i$ and $j$ such that $a_i, b_i, a_j, b_j$ are not pairwise distinct. Then, since every vertex is incident to exactly one red edge and the edges $a_ib_i, a_jb_j$ red edges, it follows that the set $\{a_i,b_i,a_j,b_j\}$ contains exactly two distinct vertices and in particular, that $a_i \in \{a_j, b_j\}$. Then $y_i a_i$ and $a_i y_k$, where $k$ is either $j$ or $j+2$, are two consecutive blue edges, and thus, once more, there exists a $7$-cycle containing them, along with three other red edges. This cycle contains a subpath $(x_i, y_i, a_i, y_k, x_k)$. But now there has to exist a path of length $3$ from $x_i$ to $x_k$ that contains one red edge, which is impossible. Therefore all the vertices $a_i$ and $b_i$ are indeed distinct, which proves the claim. So far we have shown that $\Gamma$ contains a subgraph as in \cref{figure:An_partial2}.

\begin{figure}[h]
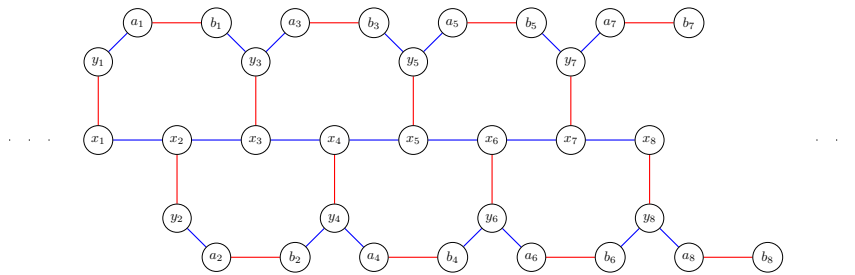

  \centering
\begin{resizedtikzpicture}{0.7\textwidth}
	\begin{pgfonlayer}{nodelayer}
		\node [style=krog] (3) at (-4.5, 0.5) {$x_1$};
		\node [style=krog] (4) at (-2.5, 0.5) {$x_2$};
		\node [style=krog] (11) at (-4.5, 2.5) {$y_1$};
		\node [style=krog] (12) at (-2.5, -1.5) {$y_2$};
		\node [style=krog] (19) at (-3.5, 3.5) {$a_1$};
		\node [style=krog] (20) at (-1.5, 3.5) {$b_1$};
		\node [style=krog] (39) at (-0.5, 0.5) {$x_3$};
		\node [style=krog] (40) at (1.5, 0.5) {$x_4$};
		\node [style=krog] (42) at (-0.5, 2.5) {$y_3$};
		\node [style=krog] (43) at (1.5, -1.5) {$y_4$};
		\node [style=krog] (45) at (0.5, 3.5) {$a_3$};
		\node [style=krog] (46) at (2.5, 3.5) {$b_3$};
		\node [style=krog] (47) at (-1.5, -2.5) {$a_2$};
		\node [style=krog] (48) at (0.5, -2.5) {$b_2$};
		\node [style=krog] (51) at (3.5, 0.5) {$x_5$};
		\node [style=krog] (52) at (5.5, 0.5) {$x_6$};
		\node [style=krog] (54) at (3.5, 2.5) {$y_5$};
		\node [style=krog] (55) at (5.5, -1.5) {$y_6$};
		\node [style=krog] (57) at (4.5, 3.5) {$a_5$};
		\node [style=krog] (58) at (6.5, 3.5) {$b_5$};
		\node [style=krog] (59) at (2.5, -2.5) {$a_4$};
		\node [style=krog] (60) at (4.5, -2.5) {$b_4$};
		\node [style=krog] (63) at (7.5, 0.5) {$x_7$};
		\node [style=krog] (64) at (9.5, 0.5) {$x_8$};
		\node [style=krog] (66) at (7.5, 2.5) {$y_7$};
		\node [style=krog] (67) at (9.5, -1.5) {$y_8$};
		\node [style=krog] (69) at (8.5, 3.5) {$a_7$};
		\node [style=krog] (70) at (10.5, 3.5) {$b_7$};
		\node [style=krog] (71) at (6.5, -2.5) {$a_6$};
		\node [style=krog] (72) at (8.5, -2.5) {$b_6$};
		\node [style=krog] (73) at (10.5, -2.5) {$a_8$};
		\node [style=krog] (74) at (12.5, -2.5) {$b_8$};
		\node [style=none] (102) at (-5, 3) {};
		\node [style=none] (103) at (-5.25, 0.5) {};
		\node [style=none] (104) at (-4.5, -1) {};
		\node [style=none] (105) at (-3, -2) {};
		\node [style=pika] (120) at (13.75, 0.5) {};
		\node [style=pika] (121) at (14.25, 0.5) {};
		\node [style=pika] (122) at (14.75, 0.5) {};
            \node [style=pika] (120) at (-5.75, 0.5) {};
		\node [style=pika] (121) at (-6.25, 0.5) {};
		\node [style=pika] (122) at (-6.75, 0.5) {};
	\end{pgfonlayer}
	\begin{pgfonlayer}{edgelayer}
		\draw [blue, thick] (11) to (19);
		\draw [red, thick](19) to (20);
		\draw [blue, thick] (20) to (42);
		\draw [blue, thick] (39) to (4);
		\draw [blue, thick] (4) to (3);
		\draw (11)[red, thick] to (3);
		\draw (12)[red, thick] to (4);
		\draw [blue, thick] (12) to (47);
		\draw[red, thick] (47) to (48);
		\draw [blue, thick] (42) to (45);
		\draw [red, thick](45) to (46);
		\draw [blue, thick](46) to (54);
		\draw [blue, thick](51) to (40);
		\draw [blue, thick](40) to (39);
		\draw [red, thick](42) to (39);
		\draw [blue, thick](48) to (43);
		\draw [red, thick](43) to (40);
		\draw [blue, thick](43) to (59);
		\draw [red, thick](59) to (60);
		\draw [blue, thick](54) to (57);
		\draw [red, thick](57) to (58);
		\draw [blue, thick](58) to (66);
		\draw [blue, thick](63) to (52);
		\draw [blue, thick](52) to (51);
		\draw [red, thick](54) to (51);
		\draw [blue, thick](60) to (55);
		\draw [red, thick](55) to (52);
		\draw [blue, thick](55) to (71);
		\draw [red, thick](71) to (72);
		\draw [blue, thick](66) to (69);
		\draw [red, thick](69) to (70);
		\draw [blue, thick](64) to (63);
		\draw [red, thick](66) to (63);
		\draw [blue, thick](72) to (67);
		\draw [red, thick](67) to (64);
		\draw [blue, thick](67) to (73);
		\draw [red, thick](73) to (74);
	\end{pgfonlayer}
\end{resizedtikzpicture}

\caption[\ ]{Vertices $x_i, y_i, a_i$ and $b_i$ in $\Gamma$ for $i \in \ZZ_n$. }
\label{figure:An_partial2}
\end{figure}

Every vertex $a_i$ must have an additional neighbour $v_i$, connected through a blue edge. We will show that $v_i =  b_{i+1}$. Since the edges $v_ia_i$ and $a_iy_i$ are both blue, the path $(v_i,a_i,y_i)$ must be contained in a girth cycle $\gamma = (v_i,a_i,y_i,w_1,w_2,w_3,w_4,v_i)$, which must contain three red edges. Then, the path of length five $P = (y_i,w_1,w_2,w_3,w_4,v_i)$ is a red-blue-red-blue-red path. By inspecting \cref{figure:An_partial2}, one can see that $w_4 \in \{b_{i-3}, a_{i-1}, b_{i-1}, a_{i+1}\}$. Since $v_i$ is connected to $w_4$ by a red edge, then necessarily $v_i \ \in \{a_{i-3}, b_{i-1}, a_{i-1}, b_{i+1}\}$. If $v_i$ is either $b_{i-1}$ or $a_{i-1}$, we get a cycle of length $6$, contradicting that the girth of $\Gamma$ is $7$. If $v_i$ is $a_{i-3}$ then, since the edges $a_iy_{i-3}$ and $y_{i-3}a_{i-3}$ are blue and consecutive, there must exist a red-blue-red-blue-red path starting at $y_{i-3}$ and ending at $a_{i}$. However, one can easily see that any such path starting at $y_{i-3}$ ends at one of $\{a_{i-4}, a_{i-6}, b_{i-4}, b_{i-1}\}$. Since $n \ge 4$, it follows that $a_i$ cannot be an element of this set. Therefore, $v_i$ must be $b_{i+1}$. By adding all edges of the form $a_ib_{i+1}$ to the graph in \cref{figure:An_partial2}, we obtain a cubic connected subgraph of $\Gamma$, which must necessarily be the entire $\Gamma$. It is straightforward to see that $\Gamma$ is isomorphic to a graph $\A(n)$ where $n = \V(\Gamma)/4$. 
\end{proof}

\begin{cor}\label{cor:lemma446}
    Let $\Gamma$ be a 
    connected cubic graph. Then $\Gamma$ is a  vertex-transitive graph of girth $7$ with $n$ vertices and signature $(4,4,6)$ if and only if $n$ is divisible by $12$ and $\Gamma$ is isomorphic to $\A(n/4)$.
\end{cor}

\begin{proof}
    Let $\Gamma$ be a connected cubic vertex-transitive graph of girth $7$ with signature $(4,4,6)$. Then by \cref{lem:446_2} it is isomorphic to $\A(m)$ where $m = n/4$. It is easy to see that for $m  \leq 7$ the graph $\A(n)$ does not have girth $7$ and signature $(4,4,6)$, and thus $m \geq 8$. Then, by \cref{lem:446_1}, $m$ must be divisible by $3$.

    Conversely, let $\Gamma$ be isomorphic to $\A(m)$ where $m$ is divisible by $12$. Then by \cref{lem:446_1}, $\Gamma$ is a vertex-transitive graph of girth $7$ and signature $(4,4,6)$.
\end{proof}

\begin{theorem}\label{th:cay}
    A connected cubic vertex-transitive graph of girth $7$ has signature $(4,4,6)$ if and only if it is isomorphic to $\text{Cay}((\ZZ_2 \times \ZZ_2) \rtimes_\varphi \ZZ_{3i}, \{(0,0,1), (0,0,-1), (1,1,0)\})$ for some $i \geq 3$, where $\varphi: \ZZ_{3i} \rightarrow \text{Aut}(\ZZ_2 \times \ZZ_2) $ maps the generator $1$ of $\ZZ_{3i}$ to the automorphism of $\ZZ_2^2$ whose matrix form in the standard basis is 
    \[\begin{bmatrix}
    1 & 1 \\
    1 & 0 \\
\end{bmatrix}.\] 
\end{theorem}

\begin{proof}
    Let us draw a subgraph of the graph $\Gamma$ in the statement, that is induced by the vertex $((0,0),0)$ and all of the vertices at distance at most $3$ from it. 

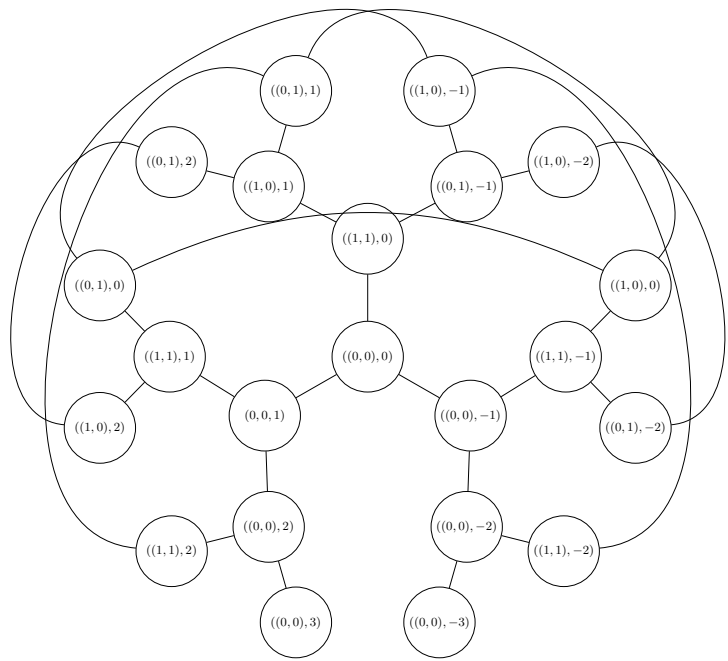
\begin{figure}[h]
\begin{adjustbox}{width=0.7\textwidth,center}
\begin{tikzpicture}[
  krog/.style={
    circle, draw, minimum size=1.8cm, inner sep=0pt, font=\small, align=center
  }, 
  povezava/.style={-}
]

\begin{scope}[rotate=-120]

  \node[krog] (0) at (0,0) {${((0,0), 0)}$};

  \node[krog] (1) at (90:3cm) {${((0,0), -1)}$};
  \node[krog] (2) at (210:3cm) {${((1,1), 0)}$};
  \node[krog] (3) at (330:3cm) {$(0,0,1)$};

  \draw[povezava] (0) -- (1);
  \draw[povezava] (0) -- (2);
  \draw[povezava] (0) -- (3);

  \node[krog] (4) at (60:5cm) {${((0,0), -2)}$};
  \node[krog] (5) at (120:5cm) {${((1,1), -1)}$};
  \node[krog] (6) at (180:5cm) {${((0,1), -1)}$};
  \node[krog] (7) at (240:5cm) {${((1,0), 1)}$};
  \node[krog] (8) at (300:5cm) {${((1,1), 1)}$};
  \node[krog] (9) at (0:5cm)   {${((0,0), 2)}$};

  \draw[povezava] (1) -- (4);
  \draw[povezava] (1) -- (5);

  \draw[povezava] (2) -- (6);
  \draw[povezava] (2) -- (7);

  \draw[povezava] (3) -- (8);
  \draw[povezava] (3) -- (9);

  \node[krog] (10) at (45:7cm)  {${((0,0), -3)}$};
  \node[krog] (11) at (75:7cm)  {${((1,1), -2)}$};
  \node[krog] (12) at (105:7cm) {${((0,1), -2)}$};
  \node[krog] (13) at (135:7cm) {${((1,0), 0)}$};
  \node[krog] (14) at (165:7cm) {${((1,0), -2)}$};
  \node[krog] (15) at (195:7cm) {${((1,0), -1)}$};
  \node[krog] (16) at (225:7cm) {${((0,1), 1)}$};
  \node[krog] (17) at (255:7cm) {${((0,1), 2)}$};
  \node[krog] (18) at (285:7cm) {${((0,1), 0)}$};
  \node[krog] (19) at (315:7cm) {${((1,0), 2)}$};
  \node[krog] (20) at (345:7cm) {${((1,1), 2)}$};
  \node[krog] (21) at (15:7cm)  {${((0,0), 3)}$};

  \draw[povezava] (4) -- (10);
  \draw[povezava] (4) -- (11);
  \draw[povezava] (5) -- (12);
  \draw[povezava] (5) -- (13);
  \draw[povezava] (6) -- (14);
  \draw[povezava] (6) -- (15);
  \draw[povezava] (7) -- (16);
  \draw[povezava] (7) -- (17);
  \draw[povezava] (8) -- (18);
  \draw[povezava] (8) -- (19);
  \draw[povezava] (9) -- (20);
  \draw[povezava] (9) -- (21);

  \draw[povezava, bend right =100] (11) to (15);
  \draw[povezava, bend right=100] (12) to (14);
  \draw[povezava, bend right=100] (13) to (16);
  \draw[povezava, bend left=100] (20) to (16);
  \draw[povezava, bend left=100] (19) to (17);
  \draw[povezava, bend left=100] (18) to (15);
  \draw[povezava, bend right=25] (13) to (18);

\end{scope}
\end{tikzpicture}
\end{adjustbox}

\caption{Vertices at distance at most $3$ in $\text{Cay}((\ZZ_2 \times \ZZ_2) \rtimes_\varphi \ZZ_{3i}, \{(0,0,\pm 1), (1,1,0)\})$.}
\label{fig:neigh3}
\end{figure}

   \cref{fig:neigh3} shows all the cycles of length at most $7$ that contain the vertex $((0,0),0)$. Observe that $\Gamma$ is a vertex-transitive graph of girth $7$ with $12i$ vertices and signature $(4,4,6)$. Hence, by Corollary~\ref{cor:lemma446}, we obtain precisely the unique vertex-transitive graph of girth $7$ and signature $(4,4,6)$ for every order for which such a graph exists.
\end{proof}

\begin{remark}
    The groups in \cref{th:cay} are isomorphic to $\langle a, b \mid a^2, b^{3i}, bab^{-1}ab^{-1}ab\rangle$ for $i \geq 3$.
\end{remark}

\section{Signatures $(4,4,4)$ $(4,5,5)$,  $(4,6,6)$ and $(5,5,6)$}\ \label{sec:455-466-556}
In this section we will deal with the four remaining possible signatures. Signatures $(4,6,6)$ and $(5,5,6)$ yield no vertex-transitive graphs, while the other two together yield exactly four graphs.

\begin{lemma}\label{lem:444}
    Let $\Gamma$ be a connected cubic vertex-transitive graph of girth $7$ with signature $(4,4,4)$. Then $\Gamma$ is arc-transitive.
\end{lemma}

\begin{proof}
    Suppose for a contradiction that there exist two orbits for the action of the automorphism group of $\Gamma$ on the edges. Then $\Gamma$ and an orbit with $|E(\Gamma)|/3$ edges cannot satisfy \cref{CONDITION}. Hence by \cref{cor:condition} we arrive to a contradiction, and $\Gamma$ is edge-transitive. Since every vertex- and edge-transitive graph of odd valence is arc-transitive, we conclude that $\Gamma$ is arc-transitive.
\end{proof}

Feng and Nedela classified cubic arc-transitive graphs of girth at most $7$ in \cite{Feng+Nedela+2006}. Their classification leads to the following proposition.

\begin{proposition}\label{pr:444}
    Let $\Gamma$ be a cubic vertex-transitive graph of girth $7$ and signature $(4,4,4)$. Then $\Gamma$ is the Coxeter graph.
\end{proposition}

\begin{proof}
    By \cref{lem:444} $\Gamma$ is arc-transitive.
    It was proven in \cite[Lemma 4.3]{Feng+Nedela+2006} that if $\Gamma$ is  a connected arc-transitive cubic graph of girth $7$ such that there are more than two $7$-cycles passing through an edge of a graph (a condition that is satisfied by $\Gamma$, since the signature is $(4,4,4)$), then it must be the Coxeter graph. It is straightforward to check that the Coxeter graph is indeed a cubic vertex-transitive graph of girth $7$ and signature $(4,4,4)$.
\end{proof}

\begin{lemma}\label{lem:not446}
    Let $\Gamma$ be a cubic vertex-transitive graph of girth $7$. Then the signature of $\Gamma$ is not $(4,6,6)$. 
\end{lemma}

\begin{proof}
    Suppose for a contradiction that $\Gamma$ is a cubic vertex-transitive graph of girth $7$ with signature $(4,6,6)$. Then the edges with $\epsilon(e) = 4$ form an $\Aut(\Gamma)$-orbit of edges $\O$ such that $|\O|=|\E(\Gamma)|/3$. We will again call the edges from this orbit \emph{red} edges, and the other edges \emph{blue}. The only integers that satisfy \cref{CONDITION} and \cref{cor:condition} are $r_1 = 0$ and $r_2 = \ldots = r_8 = 2$. That is, for an arbitrary vertex $v$, there exist $8$ girth-cycles $c_1, \ldots, c_8$ that contain it, one of the cycles contains no red edges, and the rest of them contain two each. It follows that $\Gamma - \O$ is a disjoint union of $7$-cycles. Let $\alpha = (v_0, \ldots, v_6, v_0)$ be a cycle in $\Gamma$ that contains two red edges. Assume first that $\alpha$ is as in \cref{figure:c2}, that is, assume that precisely the edges $v_0 v_1$ and $v_2 v_3$ are red (see \cref{fig:2cycl-a}). Then there exist vertices $x$ and $y$ such that $(v_0, x, y, v_3, v_4, v_5, v_6, v_0)$ is a $7$-cycle with only blue edges. But then $(v_0, x, y, v_3, v_2, v_1,v_0)$ is a $6$-cycle, which is a contradiction. It follows that every $7$-cycle in $\Gamma$ is as in \cref{figure:c3}. Without loss of generality we can assume that precisely the edges $v_0v_1$ and $v_3v_4$ in $\alpha$ are red.

\begin{figure}[h]
\begin{subfigure}[b]{0.45\textwidth} 
\centering
\begin{tikzpicture}[scale=0.8] 
\node[draw, circle, minimum size=0.65cm] (y) at ({3*cos(0)},{3*sin(0)}) {$y$};
\node[draw, circle, minimum size=0.65cm] (z) at ({3*cos(51.43)},{3*sin(51.43)}) {$v_3$};
\node[draw, circle, minimum size=0.65cm] (v4) at ({3*cos(102.86)},{3*sin(102.86)}) {$v_4$};
\node[draw, circle, minimum size=0.65cm] (v5) at ({3*cos(154.29)},{3*sin(154.29)}) {$v_5$};
\node[draw, circle, minimum size=0.65cm] (v6) at ({3*cos(205.71)},{3*sin(205.71)}) {$v_6$};
\node[draw, circle, minimum size=0.65cm] (v0) at ({3*cos(257.14)},{3*sin(257.14)}) {$v_0$};
\node[draw, circle, minimum size=0.65cm] (x) at ({3*cos(308.57)},{3*sin(308.57)}) {$x$};

\draw[blue, thick] (y) -- (z); 
\draw[blue, thick] (v4) -- (z); 
\draw[blue, thick] (v4) -- (v5) -- (v6) -- (v0); 
\draw[blue, thick] (x) -- (v0); 
\draw[blue, thick] (x) -- (y);

\node[draw, circle, minimum size=0.65cm] (v1) at (-0.3,-1) {$v_1$};
\node[draw, circle, minimum size=0.65cm] (v3) at (0.55,0.9) {$v_2$};

\draw[red, thick]  (v0) -- (v1);
\draw[blue, thick] (v3) -- (v1);
\draw[red, thick]  (v3) -- (z);

\end{tikzpicture}
\caption{\ }
\label{fig:2cycl-a}
\end{subfigure}
\hspace{3pt}
\begin{subfigure}[b]{0.45\textwidth} 
\centering
\begin{tikzpicture}[scale=0.8] 
\node[draw, circle, minimum size=0.65cm] (y) at ({3*cos(0)},{3*sin(0)}) {$y$};
\node[draw, circle, minimum size=0.65cm] (z) at ({3*cos(51.43)},{3*sin(51.43)}) {$z$};
\node[draw, circle, minimum size=0.65cm] (v4) at ({3*cos(102.86)},{3*sin(102.86)}) {$v_4$};
\node[draw, circle, minimum size=0.65cm] (v5) at ({3*cos(154.29)},{3*sin(154.29)}) {$v_5$};
\node[draw, circle, minimum size=0.65cm] (v6) at ({3*cos(205.71)},{3*sin(205.71)}) {$v_6$};
\node[draw, circle, minimum size=0.65cm] (v0) at ({3*cos(257.14)},{3*sin(257.14)}) {$v_0$};
\node[draw, circle, minimum size=0.65cm] (x) at ({3*cos(308.57)},{3*sin(308.57)}) {$x$};

\draw[blue, thick] (y) -- (z); 
\draw[blue, thick] (v4) -- (z); 
\draw[blue, thick] (v4) -- (v5) -- (v6) -- (v0); 
\draw[blue, thick] (x) -- (v0); 
\draw[blue, thick] (x) -- (y);

\node[draw, circle, minimum size=0.65cm] (v1) at (-0.25,-1.5) {$v_1$};
\node[draw, circle, minimum size=0.65cm] (v2) at (0,0) {$v_2$};
\node[draw, circle, minimum size=0.65cm] (v3) at (-0.25,1.5) {$v_3$};

\draw [red, thick]  (v0) -- (v1);
\draw[blue, thick] (v1) -- (v2);
\draw[blue, thick] (v2) -- (v3);
\draw [red, thick]  (v3) -- (v4);

\end{tikzpicture}
\caption{\ }
\label{fig:2cycl-b}
\end{subfigure}
\caption{Two cycles in a graph of signature $(4,6,6)$.}
\end{figure}
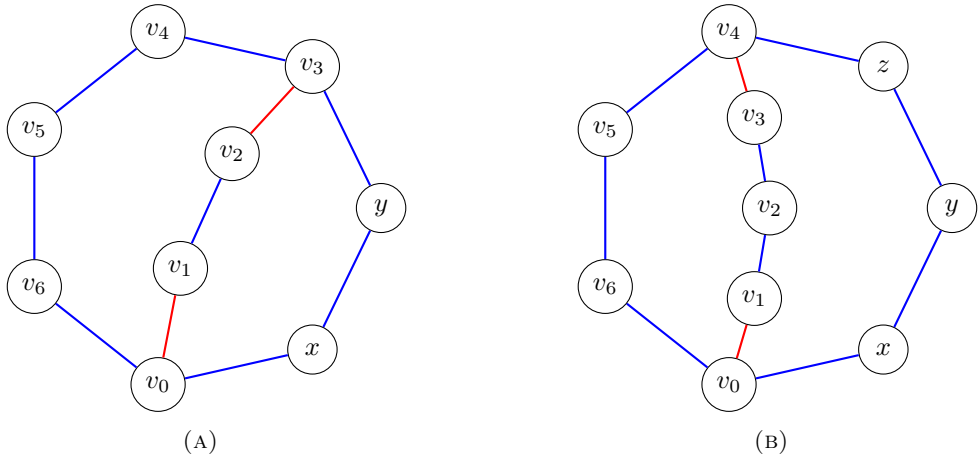

Similarly as before, there must exist vertices $x, y$ and $z$ such that $\beta = (v_0, x, y, z, v_4, v_5, v_6, v_0)$ is a $7$-cycle that contains only blue edges (see \cref{fig:2cycl-b}). 
First note that $v_2$ is not a neighbour of any vertex in $\beta$, or $\Gamma$ would contain a cycle of length strictly less than $7$.
Let $H$ be the subgroup of $\Gamma$ that preserves the cycle $\beta$. Now notice that since $\Gamma$ is vertex-transitive, and $\Aut(\Gamma)$ preserves the set of blue edges, the vertices of $\beta$ form a block of imprimitivity for the action of $\Aut(\Gamma)$ on the vertices. It follows that $H$ induces a transitive subgroup of $\Dih(\beta)$ on $\beta$. Since $7$ is odd, such a group will contain the cyclic subgroup generated by a four-step rotation, let us denote it $\varphi$, that maps $v_0$ to $v_4$. Hence $\varphi$ maps $v_1$ to $v_3$, and $\gamma = (v_1, v_2, v_1^\varphi, v_2^\varphi, \ldots, v_1^{\varphi^6}, v_2^{\varphi^6}, v_1)$ is a closed walk of length $14$ that contains only blue edges. Furthermore, no two vertices at distance $2$ in $\gamma$ can be the same (if two vertices at distance two that do not have a neighbour in $\beta$ were the same vertex, it would lie on a cycle of length $5$), as well as vertices at distance $7$ (precisely one of them is a neighbour to a vertex in $\beta$). But, since the blue edges induce a disjoint union of cycles of length $7$, we see that $\gamma$ is a subwalk of some $7$-cycle. Since any such closed walk will necessarily contain two vertices at distance two or seven in the walk that are the same, we arrive to a contradiction. It follows that a graph satisfying our conditions cannot exist.
\end{proof}

\begin{lemma}\label{lem:not556}
    Let $\Gamma$ be a cubic vertex-transitive graph of girth $7$. Then the signature of $\Gamma$ is not $(5,5,6)$.
\end{lemma}

\begin{figure}[h]
    \centering

\begin{tikzpicture}[
  scale=0.55,
  node style/.style={
    circle, draw=black, fill=white, text=black,
    minimum size=5mm, inner sep=0pt,
    opacity=1, text opacity=1, fill opacity=1
  },
  every path/.style={thick}
]

\foreach \i/\angle in {
  1/90, 2/38.57, 3/347.14, 4/295.71, 5/244.29, 6/192.86, 7/141.43
} {
  \coordinate (cx\i) at (\angle:2.6cm);
}

\foreach \i/\angle in {
  1/90, 2/38.57, 3/347.14, 4/295.71, 5/244.29, 6/192.86, 7/141.43
} {
  \coordinate (cy\i) at (\angle:4.5cm);
}

\foreach \i/\angle in {
  1/90, 2/38.57, 3/347.14, 4/295.71, 5/244.29, 6/192.86, 7/141.43
} {
  \pgfmathsetmacro{\plus}{\angle + 10}
  \pgfmathsetmacro{\minus}{\angle - 10}
  \coordinate (ca\i) at (\plus:6cm);
  \coordinate (cb\i) at (\minus:6cm);
}

\foreach \i [evaluate=\i as \j using {mod(\i,7)+1}] in {1,...,7} {
  \draw[blue] (cx\i) -- (cx\j);
}

\foreach \i in {1,...,7} {
  \draw[red] (cx\i) -- (cy\i);
}

\foreach \i in {1,...,7} {
  \draw[blue] (cy\i) -- (ca\i);
  \draw[blue] (cy\i) -- (cb\i);
}

\foreach \i in {1,...,7} {
  \pgfmathtruncatemacro{\j}{mod(\i+1,7)+1}
  \draw[red, bend left=75] (cb\i) to (ca\j);
}

\foreach \i in {1,...,7} {
  \pgfmathtruncatemacro{\k}{mod(\i+2,7)+1} 
  \draw[blue, dotted, bend left=20] (cb\i) to (ca\k);
}

\foreach \i in {1,...,7} {
  \node[node style] at (cx\i) {$x_{\i}$};
}

\foreach \i in {1,...,7} {
  \node[node style] at (cy\i) {$y_{\i}$};
}

\foreach \i in {1,...,7} {
  \node[node style] at (ca\i) {$a_{\i}$};
  \node[node style] at (cb\i) {$b_{\i}$};
}

\end{tikzpicture}

    \caption{Graph constructed in the proof of \cref{lem:not556}.}
    \label{fig:556}
\end{figure}
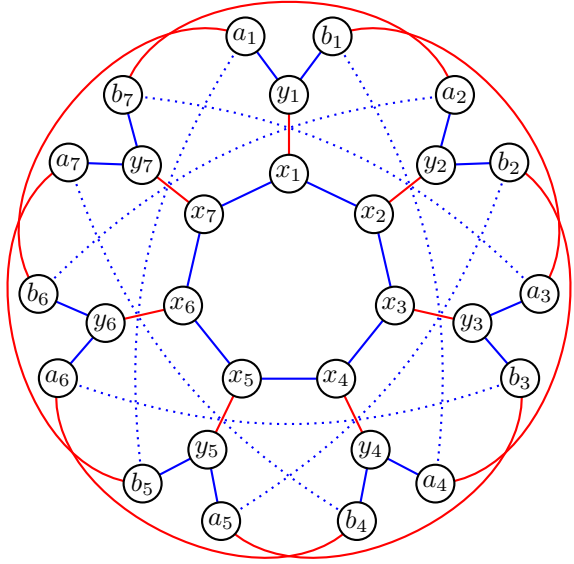

\begin{proof}
    Let $\Gamma$ be a cubic vertex-transitive graph of girth $7$ with signature $(5,5,6)$. Similarly as before, the edges $e$ with $\epsilon(e) = 6$ form an $\Aut(\Gamma)$-orbit of edges $\O$ such that every vertex in $\Gamma$ is incident to precisely one edge from $\O$. We will again call these edges \emph{red} and the rest of the edges \emph{blue}. For a $v \in V(\Gamma)$ there exist $7$ girth cycles that contain $v$.
    The only integers for which $(\Gamma, \O)$ satisfies \cref{CONDITION} are $r_2 = \ldots = r_8 = 3$, $r_1 = 0$.
     Hence, by \cref{cor:condition}, there exists $9$ $7$-cycles that contain $v$, $8$ of them contain three red edges, and one $7$-cycle $\gamma = (x_1, \ldots, x_7, x_1)$contains only blue edges. Each of the vertices $x_i$ is connected to a new vertex $y_i$ by a red edge.
    Now consider the neighbourhood of $x_1$, that is $\{x_7, x_2, y_1\}$. Let $d$ denote the number of $7$-cycles containing $(x_7, x_1, x_2)$, $e$ the number of $7$-cycles containing $(x_7, x_1, y_1)$, and $f$ the number of $7$-cycles containing $(y_1, x_1, x_2)$. We then have the system of equations
    \begin{align*}
        \begin{bmatrix}
1 & 1 & 0 \\
1 & 0 & 1 \\
0 & 1 & 1 
\end{bmatrix}  \begin{bmatrix}
d\\
e \\
f
\end{bmatrix}   = \begin{bmatrix}
5\\
5 \\
6
\end{bmatrix}.
    \end{align*}
     We can see that $e = f = 3$ and $d = 2$ is the unique solution. So there exist two $7$-cycles that contain $(x_7, x_1, x_2)$. One of them is $\gamma$. Let $\delta$ be the other one. Now, two distinct $7$ cycles ($\delta$ and $\gamma$ for instance) can have at most $3$ consecutive edges in common, for otherwise we would obtain a cycle of length $6$ or less. If $\gamma$ and $\delta$ share a $4$-path, say $(x_7, x_1, x_2,x_3)$, then $\delta$ is $(y_7,x_7, x_1, x_2,x_3,y_3,z)$ for some vertex $z$. However, the edges $zy_7$ and $y_3z$ are both blue, making $\delta$ a $7$-cycle with only two red edges, which is impossible since the number of red edges in a $7$-cycle of $\Gamma$ is either $0$ or $3$. Hence, we can assume that $\gamma$ and $\delta$ share only the path $(x_7,x_1,x_2)$ and thus there exist two adjacent vertices $b_7$ and $a_2$ such that $\delta = (y_7, x_7, x_1, x_2, y_2, b_7, a_2, y_7)$ and that the edge $b_7a_2$ is red. Note that $b_7, a_2 \notin \{x_1, \ldots, x_7, y_1, \ldots, y_7\}$, otherwise there would exist a cycle of length at most $6$. Similarly as before, there exists an automorphism $\varphi$ of $\Gamma$ such that $(x_1, \ldots, x_6, x_7)^\varphi = (x_2, \ldots, x_7, x_1)$. Let $b_i = b_7^{\varphi^i}$ and $a_{i+2} = a_2^{\varphi^i}$,  where we consider $i$ modulo $7$.

    We claim that all of the vertices $a_1, \ldots, a_7, b_1, \ldots, b_7$ are pairwise distinct. Suppose for a contradiction that they are not. Assume that $a_1 \in \{a_2, \ldots, a_7,b_1, \ldots, b_7\}$, from which we see that $a_1$ is a neighbour of two vertices $y_j$ and $y_k$ where $|j-k| \geq 3$. Then there exists a path of length $4$ between $x_j$ and $x_k$ that contains two red edges. However, there is a path of length at most $3$  with only blue edges between any two vertices in $\gamma$. That is, there exists a cycle of length $7$, containing $x_j$ and $x_k$, with precisely two red edges, which is not possible. 
    Thus the vertices $a_1, \ldots, a_7, b_1, \ldots, b_7$ are all distinct. Furthermore, since $a_2b_7$ is a red edge, we see that $a_i b_{i+5} = (a_2b_7)^{\varphi^{i-2}}$ is a red edge  for all $i \in\ZZ_7$. So far, we have shown that the graph induced by the solid (non-dotted) edges in \cref{fig:556} is a subgraph of $\Gamma$. 

     Now since $e = f = 3$, there exist three $7$-cycles that contain $(y_1, x_1, x_2)$. Since not all three of them can contain the path $(y_1,x_1,x_2,x_3)$ (for otherwise at least two of them would share $4$ consecutive edges), one of these $7$-cycles must contain the path $(y_1,x_1,x_2,y_2)$. In other words, there exists a $7$-cycle $(y_1,x_1,x_2,y_2,z_2,w,z_1, y_1)$ where $z_i \in \{a_i,b_i\}$, $w \in \V(\Gamma)$, and one of the edges $z_2w$ or $wz_1$ is red. 
     
     If $z_1 = b_1$ and $z_2=a_2$ then $w \in \{b_7,a_3\}$, since $w$ is connected to either $z_1$ or $z_2$ via a red edge. However, this implies the existence of a $6$-cycle, namely, $(b_7,y_7,x_7,x_1,y_1,b_1,b_7)$ or $(a_3,y_3,x_3,x_2,y_2,a_2,a_3)$. 

     If $z_1 = a_1$ and $z_2 = b_2$ then $w \in \{b_6,a_4\}$. Without loss of generality we can assume $w = b_6$. This implies that there is a blue edge between $b_2$ and $b_6$. However, the image of the edge $b_2b_6$ under $\varphi^{-3}$ is $b_6b_3$, implying that $b_6$ has valency $4$, a contradiction.

     Finally, if $z_1 = a_1$ and $z_2 = a_1$ (or if $z_1 = b_1$ and $z_2 = b_2$, which is an analogous case) then, since $w \neq b_7$ (for otherwise there would exist a cycle of length $6$), it follows that $w = b_6$. Then $a_2b_6$ is a blue edge whose orbit under $\langle \varphi \rangle$ consists of all the dotted edges in \cref{fig:556}. Then the graph depicted there (including the dotted edges) is the whole $\Gamma$, but $\Gamma$ is vertex-transitive by hypothesis, and the graph in \cref{fig:556} is not (some vertices are on an all-blue $7$-cycle and some are not). Therefore, this case is not possible either.

     Since we have exhausted all possibilities, the result follows.
    \end{proof}

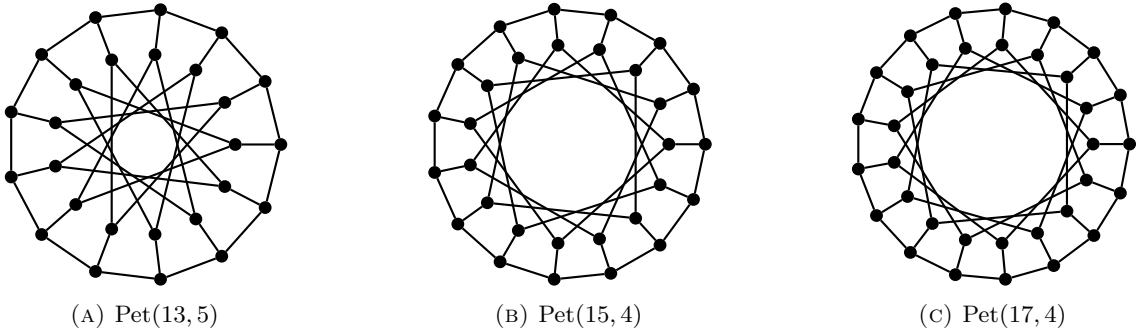
\begin{figure}[h]
            \centering
            \begin{subfigure}[b]{0.3\textwidth}
                \centering

                \begin{tikzpicture}[
                    scale=1.2,
                    every node/.style={draw, circle, minimum size=4pt, inner sep=0pt, fill=black},
                    thick
                ]
                \def\n{13}
                \def\offset{5}
                
                \def\anglestep{360/\n}
                
                \foreach \i in {0,...,\numexpr\n-1\relax} {
                    \node (o\i) at ({\anglestep*\i}:1.5) {};
                }
                
                \foreach \i in {0,...,\numexpr\n-1\relax} {
                    \node (i\i) at ({\anglestep*\i}:1) {};
                }
                
                \foreach \i in {0,...,\numexpr\n-1\relax} {
                    \pgfmathtruncatemacro{\nextv}{mod(\i+1, \n)}
                    \draw[thick] (o\i) -- (o\nextv);
                }
                
                \foreach \i in {0,...,\numexpr\n-1\relax} {
                    \draw[thick] (o\i) -- (i\i);
                }
                
                \foreach \i in {0,...,\numexpr\n-1\relax} {
                    \pgfmathtruncatemacro{\targetv}{mod(\i+\offset, \n)}
                    \draw[thick] (i\i) -- (i\targetv);
                }
            
            \end{tikzpicture}
                \caption{$\Pet(13,5)$}
                \label{fig:pet13}
            \end{subfigure}
            \hfill
            \begin{subfigure}[b]{0.3\textwidth}
                \centering
                
                \begin{tikzpicture}[
                    scale=1.2,
                    every node/.style={draw, circle, minimum size=4pt, inner sep=0pt, fill=black},
                    thick
                ]
                
                \def\n{15}
                \def\offset{4}

                \def\anglestep{360/\n}
                
                \foreach \i in {0,...,\numexpr\n-1\relax} {
                    \node (o\i) at ({\anglestep*\i}:1.5) {};
                }
                
                \foreach \i in {0,...,\numexpr\n-1\relax} {
                    \node (i\i) at ({\anglestep*\i}:1.1) {};
                }
                
                \foreach \i in {0,...,\numexpr\n-1\relax} {
                    \pgfmathtruncatemacro{\nextv}{mod(\i+1, \n)}
                    \draw[thick] (o\i) -- (o\nextv);
                }
                
                \foreach \i in {0,...,\numexpr\n-1\relax} {
                    \draw[thick] (o\i) -- (i\i);
                }
                
                \foreach \i in {0,...,\numexpr\n-1\relax} {
                    \pgfmathtruncatemacro{\targetv}{mod(\i+\offset, \n)}
                    \draw[thick] (i\i) -- (i\targetv);
                }
               
            \end{tikzpicture}
                 \caption{$\Pet(15,4)$}
                 \label{fig:pet15}
            \end{subfigure}
            \hfill
            \begin{subfigure}[b]{0.3\textwidth}
                \centering

                \begin{tikzpicture}[
                    scale=1.2,
                    every node/.style={draw, circle, minimum size=4pt, inner sep=0pt, fill=black},
                    thick
                ]
                \def\n{17}
                \def\offset{4}
                
                \def\anglestep{360/\n}
                
                \foreach \i in {0,...,\numexpr\n-1\relax} {
                    \node (o\i) at ({\anglestep*\i}:1.5) {};
                }
                
                \foreach \i in {0,...,\numexpr\n-1\relax} {
                    \node (i\i) at ({\anglestep*\i}:1.1) {};
                }
                
                \foreach \i in {0,...,\numexpr\n-1\relax} {
                    \pgfmathtruncatemacro{\nextv}{mod(\i+1, \n)}
                    \draw[thick] (o\i) -- (o\nextv);
                }
                
                \foreach \i in {0,...,\numexpr\n-1\relax} {
                    \draw[thick] (o\i) -- (i\i);
                }
                
                \foreach \i in {0,...,\numexpr\n-1\relax} {
                    \pgfmathtruncatemacro{\targetv}{mod(\i+\offset, \n)}
                    \draw[thick] (i\i) -- (i\targetv);
                } 
            
            \end{tikzpicture}
             \caption{$\Pet(17,4)$}  
             \label{fig:pet17}
            \end{subfigure}
            \caption{The three cubic vertex-transitive graphs of girth $7$ with signature $(4,5,5)$.}
            \label{fig:pet}
        \end{figure}

\begin{lemma}\label{lem:455}
    Let $\Gamma$ be a connected cubic vertex-transitive graph of girth $7$. Then the signature of $\Gamma$ is $(4,5,5)$ if and only if $\Gamma$ is isomorphic to one of the generalized Petersen graphs $\Pet(13,5)$, $\Pet(15,4)$ and $\Pet(17,4)$ shown in \cref{fig:pet}. 
\end{lemma}

\begin{proof}
    Let $\Gamma$ be a cubic vertex-transitive graph of girth $7$ with signature $(4,5,5)$. Similarly as before, the edges $e$ with $\epsilon(e) = 4$ form an $\Aut(\Gamma)$-orbit $\O$ of red edges such that every vertex in $\Gamma$ is incident to precisely one edge from $\O$. The only integers for which $(\Gamma, \O)$ satisfies \cref{CONDITION} are $r_1 = \ldots = r_7 = 2$, hence by \cref{cor:condition} each $7$-cycle contains precisely two red edges. 

    Let $v \in V(\Gamma)$, let $\Gamma(v) = \{v_1, v_2, v_3\}$ and suppose that the edge $v v_1$ is a red edge. Let $e_{1,2}$ be the number of $7$-cycles containing $(v_1, v, v_2)$, let $e_{2,3}$ be the number of $7$-cycles containing $(v_3, v, v_2)$ and let $e_{1,3}$ be the number of $7$-cycles containing $(v_1, v, v_3)$. 
    
Similarly as in the previous proof we compute that $e_{1,2} = e_{1,3} = 2$ and $e_{2,3} = 3$. Because the signature of every vertex is $(4,5,5)$, it follows that every sequence of two blue edges is contained in $3$ girth cycles, and a sequence of a blue and red edge is contained in $2$ girth cycles. 

Now let $\gamma = (x_1, \ldots, x_k, x_1)$ be a cycle that only contains blue edges. Note that $k \geq 8$ since there are no all-blue girth cycles in $\Gamma$. Since the vertices of $\gamma$ form a block of imprimitivity for $\Aut(\Gamma)$, and since $\Gamma$ is vertex-transitive, we see that that there exists an automorphism $\varphi \in \Aut(\Gamma)$ that maps $x_i$ to $x_{i+2}$, and acts on $\gamma$ as a two-step rotation. Hence if $x_i$ is connected to one of the vertices $x_1, \ldots, x_k$ with a red edge, we get a cycle of length $6$, which is a contradiction. So every vertex $x_i$ is connected to a vertex $y_i$ with a red edge, and the vertices $x_1, \ldots, x_k, y_1, \ldots, y_k$ are pairwise distinct. 

For every vertex in $\Gamma$, there exist seven $7$-cycles that contain it. Hence by Remark \ref{rem:2red} all of the cycles are as in \cref{figure:c2} or as  in \cref{figure:c3}. We separate two cases. 

\

\textbf{Case 1:}
Suppose first that all of the $7$-cycles in $\Gamma$ are as in \cref{figure:c3}, that is, the sequence of edges in every $7$-cycle is red-blue-blue-red-blue-blue-blue. The path $(x_1, x_2, x_3)$ consists of two blue edges, hence there exists three $7$-cycles containing it. Since two $7$-cycles can agree on at most $4$ consecutive vertices, and since a $7$-cycle contains at most $3$ (and never just $1$) consecutive blue edges, it follows that these three $7$-cycles contain, respectively) the paths $(y_1,x_1,x_2, x_3, y_3)$, $(y_k, x_k, x_1,x_2, x_3, y_3)$ and $(y_1,x_1,x_2, x_3, x_4, y_4)$.
In general, there exists a $7$-cycle containing any sequence of consecutive edges with colors red-blue-blue-blue-red and red-blue-blue-red. 

Hence there exists a path of length $2$ from $y_i$ to $y_{i+3}$ for any $i$, where we take $i$ modulo $k$. Let $(y_i, z_i, y_{i+3})$ be the path between $y_i$ and $y_{i+3}$. Since each $y_i$ is already incident to a red edge $x_i y_i$, and $z_i \notin \{x_1, \ldots, x_k\}$, it follows that the edges $y_i z_i$ and $z_i y_{i+3}$ are blue.  If $k \equiv 0 \pmod 3$, we get a closed walk $(y_1, z_1, y_4, z_4, \ldots, z_{k-2}, y_1)$ with at most $2k/3$ blue edges, which is a contradiction, since every cycle consisting of blue edges has length $k$. 

Therefore we see that $k \not\equiv 0 \pmod 3$. The closed walk  $(y_1, z_1, y_4, z_4, \ldots, z_{k-2}, y_1)$ that contains all the vertices $y_i$, contains at most $k$ distinct vertices, hence $z_1, \ldots, z_k \in \{y_1, \ldots, y_k\}$.
Assume that $z_1 = y_m$ and 
$z_4 = y_{\l}$ for some $m, \l$. In what follows we will examine what $m$ and $\l$ can be.
Note that $z_4 \neq z_1$, since that would imply $y_1 = y_7$, and we know that $k \geq 8$. 
Hence $(y_1, y_m, y_4, y_{\l})$ is a $3$ path, and $(x_1, y_1, y_m, y_4, y_{\l}, x_{\l})$ is a $5$ path in which precisely the first and last edge are red. Hence there exists a path consisting of two blue edges between $x_1$ and $x_{\l}$. Therefore $\l \in \{3,k-1\}$. 
If $\l = 3$, there exists an edge $y_4 y_3$, which gives a cycle of length $4$, hence we get a contradiction. It follows that $\l = k - 1$. 
Similarly $(x_{m}, y_m, y_4, y_{k-1}, x_{k-1})$ is a $4$-path in which precisely the first and the last edge belong to $\O$. Hence there exists a path consisting of three blue edges between $x_m$ and $x_{k-1}$. It follows that $m \in \{2,k-4\}$.
Similarly as before, $m \neq 2$, and thus $m = k - 4$. 

It follows that $y_4$ is connected via blue edge to both $y_{k-4}$ and $y_{k-1}$. In general, for every $j \in \{1, \ldots, k\}$, $y_j$ is connected via a blue edge to $y_{j+k-5}$, $y_{j+k+5}$, $y_{j+k-8}$ and $y_{j+k+8}$. However, each $y_j$ is connected via a blue edge to only two vertices, and thus the set of indices $\{{j+k-5}, {j+k+5}, {j+k-8},  {j+k+8}\}$ must have cardinality $2$. Since $k \geq 8$, this is only possible if $k = 13$. We have shown that $\Gamma$ is isomorphic to the generalised Petersen graph $\text{Pet}(13,5)$, shown in \cref{fig:pet13}. One can verify that it has the appropriate signature and is vertex-transitive (see also \cite[Theorem 2]{Frucht+Graver+Watkins+1971}).

\

\textbf{Case 2:} Suppose now that all of the $7$-cycles in $\Gamma$ are as in \cref{figure:c2}, that is, every $7$-cycle has edges of the sequence red-blue-red-blue-blue-blue-blue. Let $i \in \ZZ_k$. There exist three $7$-cycles that contain the walk $(x_{i-1}, x_i, x_{i+1})$. A $7$-cycle has precisely one or precisely four consecutive blue edges. Hence the  $7$-cycles that contain the walk $(x_{i-1}, x_i, x_{i+1})$ contain one of the walks $(y_{i-3}, x_{i-3}, x_{i-2}, x_{i-1}, x_i, x_{i+1}, y_{i+1})$, $(y_{i-2}, x_{i-2}, x_{i-1}, x_i, x_{i+1}, x_{i+2}, y_{i+2})$ and $(y_{i-2}, x_{i-1},x_i, x_{i+1}, x_{i+2}, x_{i+3},  y_{i+3})$. Since two distinct $7$-cycles can have at most three consecutive vertices in common, there exist precisely one $7$-cycle that contains each of the previous walks. In particular, we see that for every $i \in \ZZ_n$, the vertices $y_i$ and $y_{i+4}$ are neighbors. 
Similarly, since there exist two $7$-cycles that contain the red-blue path $(y_i, x_i, x_{i+1})$, one of them contains $(y_i, x_i, x_{i+1}, y_{i+1})$, hence for every $i \in \ZZ_k$, there exists a path of length $4$ that begins in $y_i$, ends in $y_{i+1}$, and contains only blue edges.
Recall that $k \geq 8$, since there exists no $7$-cycles that only contains blue edges. 

Suppose first that $k = 8$. Except for the edges $y_i, y_{i+4}$, there is no other edge between vertices in $\{y_1, \ldots, y_8\}$, otherwise there would exist a cycle of length at most $6$.
Hence there exist paths $(y_1, z_1, z_2, z_3, y_2)$ and $(y_1, w_1, w_2, w_3, y_8)$ with only blue edges such that $z_1, z_2, z_3, w_1, w_2, w_3 \notin \{y_1, \ldots, y_8\}$. Since $y_1$ only has one blue neighbour that is not $y_5$, it follows that $z_1 = w_1$. By the same argument $z_2 = w_2, z_3 = w_3$ and hence $x_8 = x_1$, so we arrive to a contradiction. 

Therefore $k \geq 9$. In this case every vertex is already incident to three vertices and thus $\Gamma$ must be is isomorphic to a generalised Petersen graph $\text{Pet}(k,4)$. Since there exists a path of length $4$ that begins in $y_i$ and ends in $y_{i+1}$, and contains only vertices from the set $\{y_1, \ldots, y_k \}$, it follows that $4 \cdot 4 \equiv \pm 1 \pmod k$, hence $k \in \{15, 17\}$. Both $\text{Pet}(15,4)$ and $\text{Pet}(17,4)$, depicted in Figures \ref{fig:pet15} and \ref{fig:pet17}, have the appropriate signature and are vertex-transitive (see also \cite[Theorem 2]{Frucht+Graver+Watkins+1971}).
\end{proof}

\section{Proof of ~\cref{main_theorem}}
We have examined all the cases, and can now collect them to prove \cref{main_theorem}:

\begin{proof}[Proof of \cref{main_theorem}]
    Let $\Gamma$ be a cubic vertex-transitive graph of girth $7$. Then it is girth-regular with signature $(a,b,c)$ that satisfies both \cref{P1} and, if an integer in the signature contains an integer that appears only once, the graph and the edge orbit with that signature satisfy \cref{cor:condition}. The only such signatures are $(0,1,1)$, $(2,2,2)$, $(4,4,4)$, $(4,4,6)$, $(4,5,5)$, $(4,6,6)$ and $(5,5,6)$. Graphs satisfying the last two signatures do not exist by 
    \cref{lem:not446} and \cref{lem:not556}. The characterisations of other signatures follow from \cref{P2}, \cref{P3}, \cref{pr:444}, \cref{cor:lemma446} and \cref{lem:455}.
\end{proof}

The graphs in parts $(2)$ and $(3)$ of \cref{main_theorem} are edge-girth regular, that is, every edge is contained in the same number of girth-cycles. We have proved that all of these graphs are arc-transitive. Hence, with cubic vertex-transitive graphs of girth $7$, edge-girth-regularity is a sufficient condition for arc-transitivity.

\begin{cor}
    Let $\Gamma$ be a vertex-transitive graph of girth $7$. If $\Gamma$ is edge-girth-regular, then it is arc-transitive.
\end{cor}

\section{Acknowledgements}
This paper was written as a part of the first authors doctoral thesis at the University of Ljubljana, under the supervision of Primo\v{z} Poto\v{c}nik, to whom we are grateful. We particularly appreciate his guidance and suggestions in the problem-solving strategy of this paper.

This research project was supported by Javna agencija za znanstvenoraziskovalno in inovacijsko dejavnost Republike Slovenije (ARIS), research program P1-0294 and J1-4351.

\printbibliography
\end{document}